\documentclass[12pt,twoside,reqno,openany]{amsart}

\pdfoutput=1

\usepackage{amsbsy,amscd,amsfonts,amsmath,amssymb,amsthm,color,
fancybox,fancyhdr,footmisc,graphics,graphicx,ifthen,mathrsfs,
multicol,multirow,nccrules,pdfpages,rotating,shuffle,stmaryrd,
textcomp,times,wasysym,wrapfig}

\usepackage[dvipsnames,svgnames,x11names]{xcolor}

\usepackage[all]{xy}
\usepackage[utf8]{inputenc}
\usepackage[T1]{fontenc}
\usepackage{mathtools}
\newtagform{EngelLie}[\scriptstyle]{$}{$}
\makeatletter\newcommand{\leqnomode}{\tagsleft@true}
\newcommand{\reqnomode}{\tagsleft@false}\makeatother

\newtheorem{Theorem}[equation]{Theorem}

\newtheorem{Proposition}[equation]{Proposition}

\newtheorem{Lemma}[equation]{Lemma}

\newtheorem{Corollary}[equation]{Corollary}
\newtheorem{Assertion}[equation]{Assertion}

\newtheorem{Observation}[equation]{Observation}

\theoremstyle{definition}

\newtheorem{Definition-Notation}[equation]{D\'efinition-Notation}

\newtheorem{Principle}[equation]{Principle}

\newcommand{\C}{\mathbb{C}}

\newcommand{\N}{\mathbb{N}}

\newcommand{\R}{\mathbb{R}}

\definecolor{blue}{cmyk}{1.,1.,0.,0.63}
\definecolor{red}{cmyk}{0.,1.,1.,0.63}
\definecolor{green}{cmyk}{1.,0.,1.,0.63}
\definecolor{black}{cmyk}{1.,1.,1.,1.}

\newcommand{\green}{\textcolor{green}}
\newcommand{\red}{\textcolor{red}}

\makeatletter
\renewcommand{\@fnsymbol}[1]
{\ensuremath{\ifcase#1\or $*$\or $**$\or $***$\or $****$\or $*****$
\else\@ctrerr\fi}}
\makeatother

\newcommand{\HEAD}[2]{%
\pagestyle{fancy}
\fancyhead[RO]{\tiny\sf\thepage}
\fancyhead[CO]{{\tiny\sf #1}}
\fancyhead[LE]{\tiny\sf\thepage}
\fancyhead[CE]{{\tiny\sf #2}}
\fancyfoot{}}

\numberwithin{equation}{section}

\newcommand{\Section}[1]{
\renewcommand{\thesection}{\bf\arabic{section}}
\section{#1}
\renewcommand{\thesection}{\arabic{section}}}

\newcommand{\SectionHead}[2]{
\Section{\bf #1}
\label{#2}
\HEAD{\ref{#2}.~{\sf 
#1}}{
Jo\"el {\sc Merker}, 
D\'epartement de Math\'ematiques d'Orsay, 
Universit\'e Paris-Saclay, France}}

\newcommand{\style}[1]{{\sf #1}}

\newcommand{\eqFG}{\style{eqFG}}

\newcommand{\eqL}{\style{eqL}}

\newcommand{\GL}{\style{GL}}

\renewcommand{\lim}{\style{lim}}

\newcommand{\order}{\style{order}}

\newcommand{\rank}{\style{rank}}

\newcommand{\Span}{\style{Span}}

\newcommand{\centersmallbullet}{{}_{{}^{{}^{
\scriptscriptstyle{\bullet\!}}}}}

\newcommand{\Hall}{\Hall}

\newcommand{\smallbullet}{{\scriptscriptstyle{\bullet}}}

\newcommand{\vf}{\vfill\end{document}}

\makeatletter
\def\hlinethick#1{%
\noalign{\ifnum0=`}\fi\hrule \@height #1 \futurelet
\reserved@a\@xhline}
\makeatother

\newlength{\myskip}
\begin{document}

\setcounter{section}{0}

\bigskip\bigskip


\begin{center}

{\large\bf Classification of Hessian Rank 1 Affinely Homogeneous}
\label{2-3-4-Hessian-rank-1}

\medskip

{\large\bf Hypersurfaces $H^n \subset \R^{n+1}$ 
in Dimensions $n = 2, 3, 4$}

\bigskip\bigskip

Jo\"el~{\sc Merker}\footnotemark[2]

\end{center}\bigskip

\footnotetext[1]{\,
This research was supported
in part by the Polish National Science Centre (NCN) 
via the grant number 2018/29/B/ST1/02583,
and by the Norwegian Financial Mechanism
2014--2021 via the project registration number 2019/34/H/ST1/00636.}

\footnotetext[2]{\,\,
D\'epartement de Math\'ematiques d'Orsay,
CNRS, Universit\'e Paris-Saclay, 91405 Orsay Cedex,
France, {\bf joel.merker@universite-paris-saclay.fr}}

\begin{center}
\begin{minipage}[t]{12.5cm}
\parindent 0.53cm
\footnotesize
\noindent
{\sc Abstract}.
In a previous memoir, written on the occasion of the Workshop "Complex
Analysis and Geometry", Ufa 16-19 November 2021, we showed that in
every dimension $n \geqslant 5$,  
there exists\,\,---\,\,unexpectedly\,\,---\,\,no affinely homogeneous
hypersurface $H^n \subset \R^{n+1}$ having Hessian of constant rank 1
(and not being affinely equivalent to a product with 
$\R^{m \geqslant 1}$).

The present article is devoted to determine
all non-product constant Hessian rank 1
affinely homogeneous hypersurfaces 
$H^n \subset \R^{n+1}$ in dimensions $n = 2, 3, 4$,
the cases $n = 1, 2$ being known.
Some statements of the mentioned general-dimensional memoir
are used here.

With complete details in the case $n = 2$, we illustrate the main
features of what can be termed the {\sl power series method of
equivalence}. The gist is to capture invariants at the origin
only, to create branches, and to infinitesimalize 
calculations.

In dimension $n = 3$, we find a single homogeneous model:
\[
u
\,=\,
\frac{1}{3\,z^2}
\Big\{
\big(
1-2\,y+y^2-2\,xz
\big)^{3/2}
-
(1-y)\,
\big(
1-2\,y+y^2-3\,xz
\big)
\Big\},
\]
the singularity $\frac{1}{3\,z^2}$ being illusory.

In dimension $n = 4$, 
without reaching closed forms,
we find two\,\,---\,\,depending just on some sign choice 
$\pm$\,\,---\,\,simply homogeneous models,
with their power series up to order 8,
which is sufficient to get 4 explicit affine vector fields.
\end{minipage}
\end{center}

\SectionHead{Introduction}
{introduction-2-3-4-hessian-rank-1}

The goal of this article is to determine all affinely homogeneous
local hypersurfaces $H^n \subset \R^{n+1}$ in dimensions
$n = 2, 3, 4$, the cases $n = 1, 2$ being known in the 
literature~{\cite
{Abdalla-Dillen-Vrancken-1997,
Doubrov-Komrakov-Rabinovich-1996,
Doubrov-Komrakov-1998,
Eastwood-Ezhov-1999,
Wermann-2001,
Eastwood-Ezhov-2001-2,
Olver-2018,
Merker-2019,
Chen-Merker-2019,
Arnaldsson-Valiquette-2020,
Chen-Merker-2020,
Foo-Merker-Nurowski-Ta-2021}}.
Considerations, methods, results, are also valid over $\C$.

As in~{\cite{Merker-2022}}, we graph such hypersurfaces as:
\[
u
\,=\,
F\big(x_1,x_2,x_3,x_4,\dots,x_n\big),
\]
with $F$ expandable at the origin in convergent power series.
The main hypothesis is that the $n \times n$ 
Hessian matrix
$\big( F_{x_ix_j}\big)$ has constant rank 1, 
an affinely invariant assumption~{\cite[Sec.~2]{Merker-2022}}.

In a previous article~{\cite{Chen-Merker-2020}},
handling only power-series,
Chen-Merker computed explicitly some of the 
differential invariants which naturally appear 
during the branching process, 
in order to explicitly write down the so-called
{\sl Lie-Fels-Olver recurrence relations} within each branch.
A bit before, a similar work was done 
by Arnaldsson-Valiquette~{\cite{Arnaldsson-Valiquette-2020}},
handling differential forms.

Then as a special case, Chen-Merker assumed the differential invariants
to be all constant, they examined the appearing algebraic equations,
and they re-obtained with yet another approach 
the known classification of 
affinely homogeneous nondegenerate surfaces $S^2 \subset \C^3$, 
due\,\,---\,\,more generally over $\R$\,\,---\,\,to
Abdalla-Dillen-Vrancken~{\cite{Arnaldsson-Valiquette-2020}},
Doubrov-Komrakov-Rabinovich~{\cite{Doubrov-Komrakov-Rabinovich-1996}},
Eastwood-Ezhov~{\cite{Eastwood-Ezhov-1999}}.

Such an approach through explorations of algebras of
differential invariants is in principle the most general one,
because it embraces all possible hypersurfaces, 
{\em the majority of which are not homogeneous}.
However, in higher dimensions $n \geqslant 3$, it is
delicate to handle\,\,---\,\,often
unwieldy\,\,---\,\,explicit differential
invariants.

Therefore, in this article, we employ 
a more direct and economic approach, 
which is focused only on the determination
of homogeneous models, hence disregards the complexity of
non-homogeneous geometric structures with their infinitely 
numerous differential invariants.

Presenting complete details in the case $n = 2$, we illustrate the main
features of what can be termed the {\sl power series method of
equivalence}. The gist is to capture invariants at the origin
only, to create branches, and to infinitesimalize 
calculations~{\cite{Lie-Merker-2015}}.

In dimension $n = 2$, the branching tree is the following:
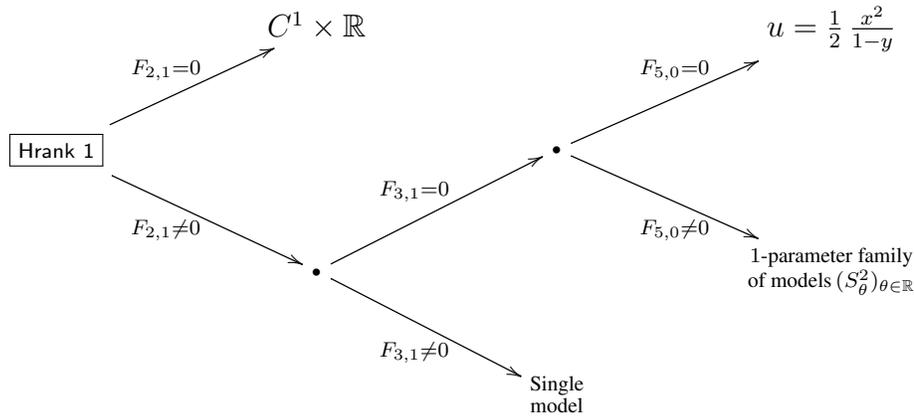
\[
\xymatrix{
&&
C^1\times\R
&&
&&
u=\tfrac{1}{2}\,\tfrac{x^2}{1-y}
\\
\ar[urr]^{F_{2,1}=0}
\boxed{\scriptstyle\sf Hrank~1}
\ar[drr]_{F_{2,1}\neq0}
&&
&&
\ar[urr]^{F_{5,0}=0}
\centersmallbullet
\ar[drr]_{F_{5,0}\neq0}
&&
\\
&&
\ar[urr]^{F_{3,1}=0}
\centersmallbullet
\ar[drr]_{F_{3,1}\neq 0}
&&
&&
{\substack{
1\text{-parameter family}
\\
\text{of models}\,(S_\theta^2)_{\theta\in\R}}}
\\
&&
&&
{\substack{
\text{Single}
\\
\text{model}}}
}
\]
We refer to Section~{\ref{surfaces-S2-R3}} for precise statements,
especially for the affine Lie algebras of the concerned (known)
homogeneous models.

In dimension $n = 3$, {\em see} 
Section~{\ref{threefolds-H3-R4}},
the 
three\footnote{\,Tree times tree! Misprints for fun!} 
is:
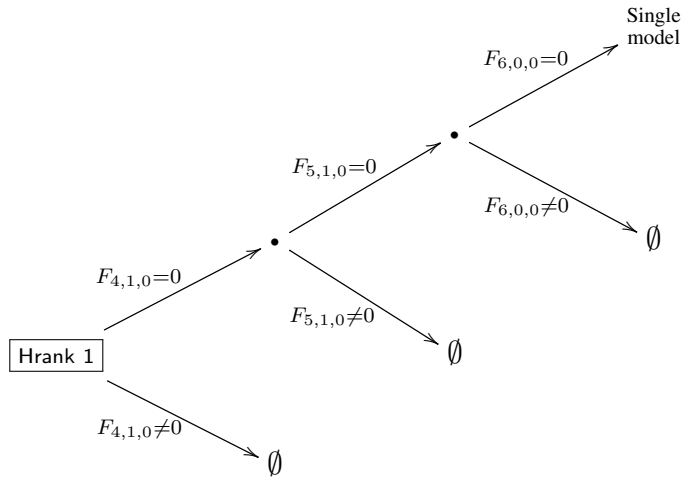
\[
\xymatrix{
&&
&&
&&
{\substack{
\text{Single}
\\
\text{model}}}
\\
&&
&&
\ar[urr]^{F_{6,0,0}=0}
\centersmallbullet
\ar[drr]_{F_{6,0,0}\neq0}
\\
&&
\ar[urr]^{F_{5,1,0}=0}
\centersmallbullet
\ar[drr]_{F_{5,1,0}\neq0}
&&
&&
\emptyset
\\
\ar[urr]^{F_{4,1,0}=0}
\boxed{\scriptstyle\sf Hrank~1}
\ar[drr]_{F_{4,1,0}\neq0}
&&
&&
\emptyset
\\
&&
\emptyset
&&
}
\]
and we find a single homogeneous model:
\[
u
\,=\,
\frac{1}{3\,z^2}
\Big\{
\big(
1-2\,y+y^2-2\,xz
\big)^{3/2}
-
(1-y)\,
\big(
1-2\,y+y^2-3\,xz
\big)
\Big\},
\]
the singularity $\frac{1}{3\,z^2}$ being illusory, with graphed 
equation:
\[
\footnotesize
\aligned
u
&
\,=\,
\ \
\tfrac{x^2}{2}
\\
&
\ \ \ \ \
+
\tfrac{x^2y}{2}
\\
&
\ \ \ \ \
+
\tfrac{x^3z}{6}
+
\tfrac{x^2y^2}{2}
\\
&
\ \ \ \ \
+
\tfrac{x^3yz}{2}
+
\tfrac{x^2y^3}{2}
\\
&
\ \ \ \ \
+
\tfrac{1}{8}\,x^4z^2
+
x^3y^2z
+
\tfrac{1}{2}\,x^2y^4
\\
&
\ \ \ \ \
+
\tfrac{5}{8}\,x^4yz^2
+
\tfrac{5}{3}\,x^3y^3z
+
\tfrac{1}{2}\,x^2y^5,
\\
&
\ \ \ \ \
+
\tfrac{1}{8}\,x^5z^3
+
\tfrac{15}{8}\,x^4y^2z^2
+
\tfrac{5}{2}\,x^3y^4z
+
\tfrac{1}{2}\,x^2y^6
\\
&
\ \ \ \ \
+
\tfrac{7}{8}\,x^5yz^3
+
\tfrac{35}{8}\,x^4y^3z^2
+
\tfrac{7}{2}\,x^3y^5z
+
\tfrac{1}{2}\,x^2y^7
\\
&
\ \ \ \ \
+
\tfrac{7}{48}\,x^6z^4
+
\tfrac{7}{2}\,x^5y^2z^3
+
\tfrac{35}{4}\,x^4y^4z^2
+
\tfrac{14}{3}\,x^3y^6z
+
\tfrac{1}{2}\,x^2y^8
+
\\
&
\ \ \ \ \
+
{\rm O}_{x,y,z}(11),
\endaligned
\]
and with affine Lie algebra:
\[
\aligned
e_1
&
\,:=\,
(1-y)\,\partial_x
-
z\,\partial_y
+
x\,\partial_u,
\\
e_2
&
\,:=\,
(1-y)\,\partial_y
-
2z\,\partial_z
+
u\,\partial_u,
\\
e_3
&
\,:=\,
u\,\partial_x
-
\tfrac{4}{3}\,x\,\partial_y
+
(1-y)\,\partial_z,
\\
e_4
&
\,:=\,
x\,\partial_x
-
z\,\partial_z
+
2\,u\,\partial_u.
\endaligned
\]

In dimension $n = 4$, 
without reaching closed forms,
we find two\,\,---\,\,depending just on some sign choice 
$\pm$\,\,---\,\,simply homogeneous models,
with their power series up to order 8,
which is sufficient to get 4 explicit affine vector fields:
\[
\aligned
e_1
&
\,:=\,
\big(
1-y\pm\tfrac{1}{5}u
\big)\,\partial_x
+
\big(
\mp\tfrac{1}{5}x-z
\big)\,\partial_y
+
\big(
-w-\tfrac{4}{75}u
\big)\,\partial_z
+
\big(
\tfrac{8}{75}x\pm\tfrac{2}{5}z
\big)\,\partial_w
+
x\,\partial_u,
\\
e_2
&
\,:=\,
-\,x\partial_w
+
(1-y)\,\partial_y
-
z\,\partial_z
-
w\,\partial_w
-
u\,\partial_u,
\\
e_3
&
\,:=\,
\tfrac{2}{3}\,u\,\partial_x
-
x\,\partial_y
+
\big(
1-y\mp\tfrac{1}{15}u
\big)\,\partial_z
+
\big(
\pm\tfrac{2}{15}x-\tfrac{2}{3}z
\big)\,\partial_w,
\\
e_4
&
\,:=\,
\pm\,\tfrac{5}{4}\,x\,\partial_x
+
\tfrac{1}{2}\,u\,\partial_y
+
\big(
-x+\tfrac{5}{4}\,z
\big)\,\partial_z
+
\big(
1-y\mp\tfrac{5}{2}w\mp\tfrac{1}{15}u
\big)\,\partial_w
\pm
\tfrac{5}{2}\,u\,\partial_u,
\endaligned
\]

The power series method of equivalence 
can be applied to other geometric structures,
including equivalences under infinite-dimensional group actions,
{\em cf.} for instance~{\cite{Eastwood-Ezhov-2004,
Eastwood-Ezhov-Isaev-2004,
Eastwood-2005,
Foo-Merker-Ta-2020,
Foo-Merker-Nurowski-Ta-2021}}.

\SectionHead{Surfaces $S^2 \subset \R^3$}
{surfaces-S2-R3}

After translation, an affine transformation of $\R^3$
fixes the origin. Consider therefore a linear map
$(x,y,u) \longmapsto (r,s,v)$:
\[
\aligned
r
&
\,:=\,
a_{1,1}\,x+a_{1,2}\,y+b_1\,u,
\\
s
&
\,:=\,
a_{2,1}\,x+a_{2,2}\,y+b_2\,u,
\\
v
&
\,:=\,
c_1\,x+c_2\,y+d\,u,
\endaligned
\ \ \ \ \ \ \ \ \ \ \ \ \ \ \ \ \ \ \ \
\text{with}
\ \ \ \ \ \ \ \ \ \ \ \ \ \ \ \ \ \ \ \
0
\,\neq\,
\left\vert\!
\begin{array}{ccc}
a_{1,1} & a_{1,2} & b_1
\\
a_{2,1} & b_{2,2} & b_2
\\
c_1 & c_2 & d
\end{array}
\!\right\vert.
\]

Also, 
consider two local analytic surfaces passing through the origin,
graphed as:
\[
u
\,=\,
F(x,y)
\ \ \ \ \ \ \ \ 
{\scriptstyle{(F(0,0)\,=\,0)}}
\ \ \ \ \ \ \ \ \ \ \ \
\text{and}
\ \ \ \ \ \ \ \ \ \ \ \ 
v
\,=\,
G(r,s)
\ \ \ \ \ \ \ \ 
{\scriptstyle{(0\,=\,G(0,0))}},
\]
with convergent series:
\[
F
\,=\,
\sum_{i+j\geqslant 1}\,
F_{j,k}\,
\frac{x^i}{i!}\,\frac{y^j}{j!}
\ \ \ \ \ \ \ \ \ \ \ \ \ \ \ \ \ \ \ \
\text{and}
\ \ \ \ \ \ \ \ \ \ \ \ \ \ \ \ \ \ \ \
G
\,=\,
\sum_{k+l\geqslant 1}\,
G_{k,l}\,
\frac{r^k}{k!}\,
\frac{s^l}{l!}.
\]

The linear map above sends the left 
surface $\{u = F\}$ to the right surface $\{v = G\}$ if and
only if the {\sl fundamental equation:}
\leqnomode\usetagform{default}
\begin{align}
\label{eqFG-S2-R3}
0
&
\,\equiv\,
\eqFG(x,y),
\end{align}
holds identically in $\R\{x,y\}$, where:
\[
\eqFG
\,:=\,
-\,c_1\,x-c_2\,y-d\,F(x,y)
+
G
\Big(
a_{1,1}x+a_{1,2}y+b_1F(x,y),\,\,
a_{2,1}x+a_{2,2}y+b_2F(x,y)
\Big),
\]
so that:
\[
0
\,=\,
\eqFG
\,=\,
\sum_{i,j\in\N}\,
\mathcal{C}_{i,j}
\Big(
a_{\smallbullet,\smallbullet},\,
b_{\smallbullet},\,
c_{\smallbullet},\,
d_{\smallbullet},\,\,
F_{\smallbullet,\smallbullet},\,
G_{\smallbullet,\smallbullet}
\Big)\,
x^i\,y^j.
\]

The {\sl core work} is to {\em compute} these
(often complicated) coefficients $\mathcal{C}_{i,j} = 0$
and to {\em analyze} their vanishing.

For all $i,j \in \N$, the coefficient of $x^i y^j$ in
$\eqFG$ can, in a standard way, be denoted as:
\[
\big[x^i\,y^j\big]
\eqFG
\,:=\,
\mathcal{C}_{i,j}
\,=\,
0,
\]
and we will constantly indicate the corresponding
indices 
$\green{\bf i, j}$ over the equal sign as:
\[
0
\overset{\green{\bf i,j}}{\,\,=\,\,}
\mathcal{C}_{i,j}.
\]
We will proceed inductively, {\sl order by order}, where:
\[
\order
\,:=\,
i+j.
\]

Two obvious affine transformations make horizontal the 
two tangent spaces:
\[
u
\,=\,
\red{\bf 0}
+
{\rm O}_{x,y}(2)
\ \ \ \ \ \ \ \ \ \ \ \ \ \ \ \ \ \ \ \
\text{and}
\ \ \ \ \ \ \ \ \ \ \ \ \ \ \ \ \ \ \ \
v
\,=\,
\red{\bf 0}
+
{\rm O}_{r,s}(2),
\]
where the two $\red{\bf 0}$ 
should be interpreted as {\sl normalizations} of order 1 terms.

\begin{Lemma}
Stabilization of order 1 terms holds if and only if $0 = c_1 = c_2$:
\[
\left[
\begin{array}{ccc}
a_{1,1} & a_{1,2} & b_1
\\
a_{2,1} & a_{2,2} & b_2
\\
c_1 & c_2 & d
\end{array}
\right]^{\green{\bf 0}}
\,\,\,\leadsto\,\,\,
\left[
\begin{array}{ccc}
a_{1,1} & a_{1,2} & b_1
\\
a_{2,1} & a_{2,2} & b_2
\\
\red{\bf 0} & \red{\bf 0} & d
\end{array}
\right]^{\green{\bf 1}}.
\]
\end{Lemma}

\proof
Apply~{\cite[Sec.~2]{Merker-2022}}, or read from~{\eqref{eqFG-S2-R3}}:
\begin{align}
0
&
\overset{\green{\bf 1,0}}{\,\,=\,\,}
-\,c_1,
\notag
\\
0
&
\overset{\green{\bf 0,1}}{\,\,=\,\,}
-\,c_2.
\qedhere
\end{align}
\endproof

Next, pass to order 2.
Possibly after rotation in the $(x,y)$-space and
in the $(r,s)$-space, 
the constant Hessian rank 1 hypothesis\,\,---\,\,which
is affinely invariant~{\cite[Sec.~2]{Merker-2022}}\,\,---\,\,reads as:
\[
F_{xx}
\,\neq\,
0
\,\equiv\,
\left\vert\!
\begin{array}{cc}
F_{xx} & F_{xy}
\\
F_{yx} & F_{yy}
\end{array}
\!\right\vert
\ \ \ \ \ \ \ \ \ \ \ \ \ \ \ \ \ \ \ \
\overset{\text{known}}{\Longleftrightarrow}
\ \ \ \ \ \ \ \ \ \ \ \ \ \ \ \ \ \ \ \
G_{rr}
\,\neq\,
0
\,\equiv\,
\left\vert\!
\begin{array}{cc}
G_{rr} & G_{rs}
\\
G_{sr} & G_{ss}
\end{array}
\!\right\vert.
\]
Thus at order 2, we have to normalize two rank 1 basic quadratic forms:
\[
F_{2,0}\,\tfrac{x^2}{2}
+
F_{1,1}\,xy
+
F_{0,2}\,\tfrac{y^2}{2}
\ \ \ \ \ \ \ \ \ \ \ \ \ \ \ \ \ \ \ \
\text{and}
\ \ \ \ \ \ \ \ \ \ \ \ \ \ \ \ \ \ \ \
G_{2,0}\,\tfrac{r^2}{2}
+
G_{1,1}\,rs
+
G_{0,2}\,\tfrac{s^2}{2}.
\]

\begin{Proposition}
Two appropriate linear transformations in the $(x,y)$-space
and in the $(r,s)$-space normalize:
\[
u
\,=\,
\tfrac{x^2}{2}
+
{\rm O}_{x,y}(3)
\ \ \ \ \ \ \ \ \ \ \ \ \ \ \ \ \ \ \ \
\text{and}
\ \ \ \ \ \ \ \ \ \ \ \ \ \ \ \ \ \ \ \
v
\,=\,
\tfrac{r^2}{2}
+
{\rm O}_{r,s}(3).
\]

Furthermore, stabilization of order $2$ terms holds if and only if
$a_{1,2} = 0$ and $d = a_{1,1}^2$:
\[
\left[
\begin{array}{ccc}
a_{1,1} & a_{1,2} & b_1
\\
a_{2,1} & a_{2,2} & b_2
\\
\red{\bf 0} & \red{\bf 0} & d
\end{array}
\right]^{\green{\bf 1}}
\,\,\,\leadsto\,\,\,
\left[
\begin{array}{ccc}
a_{1,1} & \red{\bf 0} & b_1
\\
a_{2,1} & a_{2,2} & b_2
\\
\red{\bf 0} & \red{\bf 0} & \red{a_{1,1}^2}
\end{array}
\right]^{\green{\bf 2}}.
\]
\end{Proposition}

\proof
The first assertion is known. Then the second
follows by computing the three equations
$\overset{\green{\bf 2,0}}{\,=\,}$, 
$\overset{\green{\bf 1,1}}{\,=\,}$, 
$\overset{\green{\bf 0,2}}{\,=\,}$
from~{\eqref{eqFG-S2-R3}}.
\endproof

All this is in fact proved in~{\cite[Sec.~2-5]{Merker-2022}}
for constant Hessian rank 1 hypersurfaces $H^n \subset \R^{n+1}$,
in any dimension $n \geqslant 1$.

\smallskip

Next, let order 3 monomials appear:
\[
\aligned
u
&
\,=\,
\tfrac{x^2}{2}
+
F_{3,0}\,
\tfrac{x^3}{6}
+
F_{2,1}\,
\tfrac{x^2y}{2}
+
F_{1,2}\,
\tfrac{xy^2}{2}
+
F_{0,3}\,
\tfrac{y^3}{6}
+
{\rm O}_{x,y}(4),
\\
v
&
\,=\,
\tfrac{r^2}{2}
+
G_{3,0}\,
\tfrac{r^3}{6}
+
G_{2,1}\,
\tfrac{r^2s}{2}
+
G_{1,2}\,
\tfrac{rs^2}{2}
+
G_{0,3}\,
\tfrac{s^3}{6}
+
{\rm O}_{r,s}(4).
\endaligned
\]
Then Hessian rank 1 implies (exercise) $0 = F_{1,2} = F_{0,3}$ and
$G_{1,2} = G_{0,3} = 0$:
\[
u
\,=\,
\tfrac{x^2}{2}
+
F_{3,0}\,
\tfrac{x^3}{6}
+
F_{2,1}\,
\tfrac{x^2y}{2}
+
{\rm O}_{x,y}(4)
\ \ \ \ \ \ \ \
\xrightarrow[{\rule[0pt]{50pt}{0pt}}]{\text{Equivalence}}
\ \ \ \ \ \ \ \
v
\,=\,
\tfrac{r^2}{2}
+
G_{3,0}\,
\tfrac{r^3}{6}
+
G_{2,1}\,
\tfrac{r^2s}{2}
+
{\rm O}_{r,s}(4).
\]

Now starts the real work. The fundamental equation gives:
\[
\aligned
0
&
\overset{\green{\bf 3,0}}{\,\,=\,\,}
-\,a_{1,1}^2\,F_{3,0}
+
a_{1,1}^3\,G_{3,0}
+
3\,a_{1,1}^2\,a_{2,1}\,G_{2,1}
+
3\,a_{1,1}\,\boxed{b_1},
\\
0
&
\overset{\green{\bf 2,1}}{\,\,=\,\,}
-\,a_{1,1}^2\,F_{2,1}
+
a_{1,1}^2\,a_{2,2}\,G_{2,1}.
\endaligned
\]
while 
$\overset{\green{\bf 1,2}}{\,=\,}$ and  
$\overset{\green{\bf 0,3}}{\,=\,}$ 
bring nothing for they both
reduce to $0 = 0$.

Observe that since the stability group at order 2 is a subgroup
of $\GL(3,\R)$:
\[
0
\,\neq\,
\left\vert\!
\begin{array}{ccc}
a_{1,1} & 0 & b_1
\\
a_{2,1} & a_{2,2} & b_2
\\
0 & 0 & a_{1,1}^2
\end{array}
\!\right\vert
\,=\,
a_{1,1}\,a_{2,2}\,a_{1,1}^2,
\]
we have $a_{1,1} \neq 0$, and therefore, the boxed 
free group parameter $\boxed{b_1}$
can be used to normalize:
\[
G_{3,0}
\,:=\,
\red{\bf 0},
\]
just by assigning:
\[
b_1
\,:=\,
\tfrac{1}{3}\,
a_{1,1}\,
F_{3,0}
-
\red{\bf 0}
-
a_{1,1}\,a_{2,1}\,G_{2,1},
\]
replacing of course $G_{2,1} = \frac{1}{a_{2,2}}\, F_{2,1}$ from
$\overset{\green{\bf 2,1}}{\,=\,}$.

Once $G_{3,0} = 0$ is so normalized, we restart from
the surface on the right $\{v = G\}$, we place it on the left,
we change notation $(r,s,v) \longmapsto (x,y,u)$, 
$G \longmapsto F$, 
we rename it $\{ u = F\}$ thus
with $F_{3,0} = 0$,
we take another affine equivalence 
to another surface $\{ v = G\}$ on the right,
and we again normalize similarly $G_{3,0} = 0$. 

Thus without any further work,
we can assume $F_{3,0} = 0 = G_{3,0}$, 
{\em simultaneously}.

Generally, once a normalization has been made on the right,
always, 
it can also be made {\em exactly the same} on the left.

\begin{Principle}
\label{Principle-F-G}
{\sl At any order, every performed normalization will always 
be instantly achieved on both hypersurfaces
$\{u = F\}$ and $\{v = G\}$.}\qed
\end{Principle}

Thus:
\[
u
\,=\,
\tfrac{x^2}{2}
+
F_{2,1}\,
\tfrac{x^2y}{2}
+
{\rm O}_{x,y}(4)
\ \ \ \ \ \ \ \
\xrightarrow[{\rule[0pt]{50pt}{0pt}}]{\text{Equivalence}}
\ \ \ \ \ \ \ \
v
\,=\,
\tfrac{r^2}{2}
+
G_{2,1}\,
\tfrac{r^2s}{2}
+
{\rm O}_{r,s}(4).
\]
Next, since $a_{1,1} \neq 0 \neq a_{2,2}$, the remaining 
equation $\overset{\green{\bf 2,1}}{\,=\,}$, namely:
\[
\overset{\green{\bf 2,1}}{\,\,=\,\,}
-\,a_{1,1}^2\,F_{2,1}
+
a_{1,1}^2\,a_{2,2}\,G_{2,1},
\]
shows that $G_{2,1}$ is a nonzero multiple of $F_{2,1}$
This means that $F_{2,1}$ is a {\em relative invariant}.

Consequently, if we abbreviate:
\[
\boxed{\scriptstyle\sf Hrank~1}
\,:=\,
\boxed{
\aligned
0
&
\,\neq\,
F_{xx}
\\
0
&
\,\equiv\,
F_{xx}\,F_{yy}
-
F_{xy}^2
\endaligned},
\]
we must open two branches:
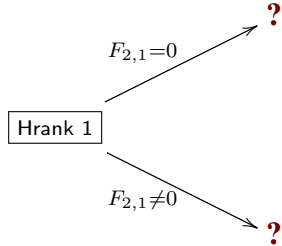
\[
\xymatrix{
&&
\red{\text{\bf ?}}
\\
\ar[urr]^{F_{2,1}=0}
\boxed{\scriptstyle\sf Hrank~1}
\ar[drr]_{F_{2,1}\neq0}
\\
&&
\red{\text{\bf ?}}
}
\]

\begin{Proposition}
\label{Prp-product}
If a surface $S^2 \subset \R^3$ is affinely homogeneous and belongs
to the branch $F_{2,1} = 0$, then $F = F(x)$ depends
only on $x$, and the surface $S^2 = C^1 \times \R_y^1$
is a cylinder over a 
curve $C^1 := \{u = F(x)\}$ 
which is affinely homogeneous in $\R^2$.
\end{Proposition}

Here and below, we will disregard such degenerate situations.
That is, we will not attempt to expressly classify 
affinely homogeneous cylinders,
because the task essentialy boils down to lower dimension.

To prove this proposition, the key argument is to infinitesimalize
and to exploit transitivity. 

A general affine vector field writes:
\[
\aligned
L
&
\,=\,
\ \ 
\big(
T_1+A_{1,1}\,x+A_{1,2}\,y+B_1\,u
\big)\,\frac{\partial}{\partial x}
\\
&
\ \ \ \ \
+
\big(
T_2+A_{2,1}\,x+A_{2,2}\,y+B_2\,u
\big)\,\frac{\partial}{\partial y}
\\
&
\ \ \ \ \
+
\big(
T_0+C_1\,x+C_2\,y+D\,u
\big)\,\frac{\partial}{\partial u}.
\endaligned
\]
It is tangent to $\{u = F(x,y)\}$ if and only if:
\[
\aligned
0
&
\,\equiv\,
\eqL(x,y)
\\
&
\,=:\,
L
\big(
-\,u+F(x,y)
\big)
\Big\vert_{u=F(x,y)},
\endaligned
\]
identically as power series in $\R\{x,y\}$.
With increasing orders $\mu = 0, 1, 2, 3, \dots$, 
this $\eqL$ may be expanded:
\[
\eqL
\,=\,
\sum_{\mu=0}^\infty\,\,
\sum_{i+j=\mu}\,
{\sf Coefficient}_{i,j}\,
x^i\,y^j.
\]

As for $\eqFG$, denote:
\[
\big[x^i\,y^j\big]
\eqL
\,:=\,
\mathcal{C}_{i,j}
\,=\,
0,
\]
or shortly:
\[
0
\overset{\green{\bf i,j}}{\,\,=\,\,}
\mathcal{C}_{i,j}.
\]

Such a vector field $L$ is tangent to: 
\[
u
\,=\,
\tfrac{x^2}{2} 
+ 
{\rm O}_{x,y}(3),
\]
if and only if:
\[
\aligned
0
&
\overset{\green{\bf 0,0}}{\,\,=\,\,}
-\,T_0,
\\
0
&
\overset{\green{\bf 1,0}}{\,\,=\,\,}
-\,C_1+T_1,
\\
0
&
\overset{\green{\bf 0,1}}{\,\,=\,\,}
-\,C_2.
\endaligned
\]
We then solve these 3 equations as:
\[
T_0
\,:=\,
0,
\ \ \ \ \ \ \ \ \ \ \ \ \ \ \ \ \ \ \ \
C_1
\,:=\,
T_1,
\ \ \ \ \ \ \ \ \ \ \ \ \ \ \ \ \ \ \ \
C_2
\,:=\,
0.
\]

In fact, the key constraint of transitivity:
\[
\Span\,
\big(
\tfrac{\partial}{\partial x},\,
\tfrac{\partial}{\partial y}
\big)
\,=\,
T_{\sf origin}S
\,=\,
\Span\,
L
\big\vert_{\sf origin}
\,=\,
\Span\,
\Big(
T_1\,
\tfrac{\partial}{\partial x}
+
T_2\,
\tfrac{\partial}{\partial y}
\Big),
\]
forces to always keep $T_1, T_2$ absolutely free\,\,---\,\,never
solved.

Next, such an $L$ is tangent to:
\[
u
\,=\,
\tfrac{x^2}{2} 
+
F_{2,1}\,\tfrac{x^2y}{2}
+
{\rm O}_{x,y}(4),
\]
if and only if moreover:
\[
\aligned
0
&
\overset{\green{\bf 2,0}}{\,\,=\,\,}
-\,\tfrac{1}{2}\,D+A_{1,1}+\tfrac{1}{2}\,F_{2,1}\,T_2,
\\
0
&
\overset{\green{\bf 1,1}}{\,\,=\,\,}
A_{1,2}+F_{2,1}\,T_1,
\\
0
&
\overset{\green{\bf 0,2}}{\,\,=\,\,}
0.
\endaligned
\]
We solve:
\[
\aligned
A_{1,2}
&
\,:=\,
-\,F_{2,1}\,T_1,
\\
D
&
\,:=\,
F_{2,1}\,T_2
+
2\,A_{1,1}.
\endaligned
\]

\proof[Proof of Proposition~{\ref{Prp-product}}]
Since $F_{2,1} = 0$ is assumed, we have by letting order 4 monomials
appear:
\[
u
\,=\,
\tfrac{x^2}{2}
+
\red{\bf 0}
+
F_{4,0}\,\tfrac{x^4}{24}
+
F_{3,1}\,\tfrac{x^3y}{6}
+
F_{2,2}\,\tfrac{x^2y^2}{4}
+
F_{1,3}\,\tfrac{xy^3}{6}
+
F_{0,4}\,\tfrac{y^4}{24}
+
{\rm O}_{x,y}(5),
\]
and:
\[
\aligned
A_{1,2}
&
\,=\,
0,
\\
D
&
\,=\,
0
+
2\,A_{1,1}.
\endaligned
\]

Then at order 3, $\eqL$ gives:
\[
\aligned
0
&
\overset{\green{\bf 3,0}}{\,\,=\,\,}
\tfrac{1}{6}\,
F_{4,0}\,T_1
+
\tfrac{1}{2}\,B_1
+
\tfrac{1}{6}\,
F_{3,1}\,T_2,
\\
0
&
\overset{\green{\bf 2,1}}{\,\,=\,\,}
\tfrac{1}{2}\,
F_{3,1}\,T_1
+
\tfrac{1}{2}\,
F_{2,2}\,T_2,
\\
0
&
\overset{\green{\bf 1,2}}{\,\,=\,\,}
\tfrac{1}{2}\,
F_{2,2}\,T_1
+
\tfrac{1}{2}\,
F_{1,3}\,T_2,
\\
0
&
\overset{\green{\bf 0,3}}{\,\,=\,\,}
\tfrac{1}{6}\,
F_{1,3}\,T_1
+
\tfrac{1}{6}\,
F_{0,4}\,T_2.
\endaligned
\]
Since there can be no linear relation between
the transitivity parameters $\{T_1, T_2\}$,
we necessarily have:
\[
0
\,=\,
F_{0,4}
\,=\,
F_{1,3}
\,=\,
F_{2,2}
\,=\,
F_{3,1}.
\]

Thus:
\[
u
\,=\,
\tfrac{x^2}{2} 
+
\red{\bf 0}
+
F_{4,0}\,\tfrac{x^4}{24}
+
F_{5,0}\,\tfrac{x^5}{120}
+
F_{4,1}\,\tfrac{x^4y}{24}
+
F_{3,2}\,\tfrac{x^3y^2}{12}
+
F_{2,3}\,\tfrac{x^2y^3}{12}
+
F_{1,4}\,\tfrac{xy^4}{24}
+
F_{0,5}\,\tfrac{y^5}{120}
+
{\rm O}_{x,y}(6).
\]

Again, $\eqL$ at order 4 gives:
\[
\aligned
0
&
\overset{\green{\bf 4,0}}{\,\,=\,\,}
\tfrac{1}{12}\,F_{4,0}\,A_{1,1}
+
\tfrac{1}{24}\,F_{5,0}\,T_1
+
\tfrac{1}{24}\,F_{4,1}\,T_2,
\\
0
&
\overset{\green{\bf 3,1}}{\,\,=\,\,}
\tfrac{1}{6}\,F_{4,1}\,T_1
+
\tfrac{1}{6}\,F_{3,2}\,T_2,
\\
0
&
\overset{\green{\bf 2,2}}{\,\,=\,\,}
\tfrac{1}{4}\,F_{3,2}\,T_1
+
\tfrac{1}{4}\,F_{2,3}\,T_2,
\\
0
&
\overset{\green{\bf 1,3}}{\,\,=\,\,}
\tfrac{1}{6}\,F_{2,3}\,T_1
+
\tfrac{1}{6}\,F_{1,4}\,T_2,
\\
0
&
\overset{\green{\bf 0,4}}{\,\,=\,\,}
\tfrac{1}{24}\,F_{1,4}\,T_1
+
\tfrac{1}{24}\,F_{0,5}\,T_2.
\endaligned
\]
By freeness of $\{T_1, T_2\}$, it is necessary that:
\[
0
\,=\,
F_{0,5}
\,=\,
F_{1,4}
\,=\,
F_{2,3}
\,=\,
F_{3,2}
\,=\,
F_{4,1}.
\]

An elementary induction on the order
$\mu \geqslant 6$ shows that in the expansion:
\[
u
\,=\,
\tfrac{x^2}{2} 
+
\red{\bf 0}
+
F_{4,0}\,\tfrac{x^4}{24}
+
F_{5,0}\,\tfrac{x^5}{120}
+
\sum_{\mu=6}^\infty\,
\sum_{i+j=\mu}\,
F_{i,j}\,
\tfrac{x^i}{i!}\,
\tfrac{y^j}{j!},
\]
all $F_{i,j}$ with $j \geqslant 1$ must be zero,
so that $F = F(x)$ is in conclusion independent of $y$.

Lastly, it can be verified
that affine homogeneity in $\R^3$ of
the cylindrical surface $\big\{ (x,y,u) \colon 
u = F(x) \}$ is equivalent to
affine homogeneity in $\R^2$ of 
the curve $\big\{ (x,u) \colon 
u = F(x) \big\}$.
\endproof

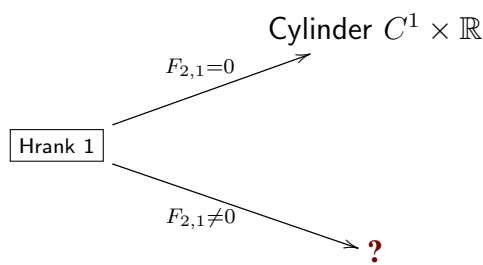
\[
\xymatrix{
&&
{\sf Cylinder}\,\,C^1\times\R
\\
\ar[urr]^{F_{2,1}=0}
\boxed{\scriptstyle\sf Hrank~1}
\ar[drr]_{F_{2,1}\neq0}
\\
&&
\red{\text{\bf ?}}
}
\]

The branch $F_{2,1} = 0$ being thus settled, assume $F_{2,1} \neq 0$.
Since $F_{2,1} \propto G_{2,1}$, 
this is a {\em coordinate-independent assumption.}
Indeed, recall:
\[
\overset{\green{\bf 2,1}}{\,\,=\,\,}
-\,a_{1,1}^2\,F_{2,1}
+
a_{1,1}^2\,a_{2,2}\,G_{2,1},
\]
with $a_{1,1} \neq 0 \neq a_{2,2}$ thanks to:
\[
0
\,\neq\,
\left\vert\!
\begin{array}{ccc}
a_{1,1} & 0 & b_1
\\
a_{2,1} & a_{2,2} & b_2
\\
0 & 0 & a_{1,1}^2
\end{array}
\!\right\vert
\,=\,
a_{1,1}\,a_{2,2}\,a_{1,1}^2.
\]

In this equation $\overset{\green{\bf 2,1}}{\,=\,}$,
it is clear that one can normalize $G_{2,1} := 1$ by choosing
$a_{2,2} := F_{2,1}$. 

In accordance with 
Principle~{\ref{Principle-F-G}}, restart, 
rename $G := F$ with $F_{2,1} = 1$, 
and normalize similarly $G_{2,1} := 1$.
Thus:
\[
u
\,=\,
\tfrac{x^2}{2}
+
\tfrac{x^2y}{2}
+
{\rm O}_{x,y}(4)
\ \ \ \ \ \ \ \
\xrightarrow[{\rule[0pt]{50pt}{0pt}}]{\text{Equivalence}}
\ \ \ \ \ \ \ \
v
\,=\,
\tfrac{r^2}{2}
+
\tfrac{r^2s}{2}
+
{\rm O}_{r,s}(4).
\]

\begin{Lemma}
Stabilization of these order $\leqslant 3$ normalizations holds 
if and only if:
\[
\left[
\begin{array}{ccc}
a_{1,1} & \red{\bf 0} & b_1
\\
a_{2,1} & a_{2,2} & b_2
\\
\red{\bf 0} & \red{\bf 0} & \red{a_{1,1}^2}
\end{array}
\right]^{\green{\bf 2}}
\,\,\,\leadsto\,\,\,
\left[
\begin{array}{ccc}
a_{1,1} & \red{\bf 0} & \red{-a_{1,1}a_{2,1}}
\\
a_{2,1} & \red{\bf 1} & b_2
\\
\red{\bf 0} & \red{\bf 0} & \red{a_{1,1}^2}
\end{array}
\right]^{\green{\bf 3}}.
\]
\end{Lemma}

\proof
Examine $\eqFG$ at order $3$:
\begin{align}
0
&
\overset{\green{\bf 3,0}}{\,\,=\,\,}
\tfrac{1}{2}\,a_{1,1}\,b_1
+
\tfrac{1}{2}\,a_{1,1}^2\,a_{2,1},
\notag
\\
0
&
\overset{\green{\bf 2,1}}{\,\,=\,\,}
-\,\tfrac{1}{2}\,
a_{1,1}^2
+
\tfrac{1}{2}\,a_{1,1}^2\,a_{2,2},
\notag
\\
0
&
\overset{\green{\bf 1,2}}{\,\,=\,\,}
0,
\notag
\\
0
&
\overset{\green{\bf 0,3}}{\,\,=\,\,}
0.
\qedhere
\end{align}
\endproof

Next, pass to order 4:
\[
\aligned
u
&
\,=\,
\tfrac{x^2}{2}
+
\tfrac{x^2y}{2}
+
F_{4,0}\,
\tfrac{x^4}{24}
+
F_{3,1}\,
\tfrac{x^3y}{6}
+
\underline{
\tfrac{x^2y^2}{2}
+
0
+
0}
+
{\rm O}_{x,y}(5),
\\
v
&
\,=\,
\tfrac{r^2}{2}
+
\tfrac{r^2s}{2}
+
G_{4,0}\,
\tfrac{r^4}{24}
+
G_{3,1}\,
\tfrac{r^3s}{6}
+
\underline{
\tfrac{r^2s^2}{2}
+
0
+
0}
+
{\rm O}_{r,s}(5).
\endaligned
\]
Here, the values of the underlined monomials are obtained 
from the (affinely invariant) hypothesis of constant
Hessian rank 1.

Indeed, from:
\[
F_{yy}
\,\equiv\,
\frac{F_{xy}^2}{F_{xx}}
\ \ \ \ \ \ \ \ \ \ \ \ \ \ \ \ \ \ \ \
\Longleftrightarrow
\ \ \ \ \ \ \ \ \ \ \ \ \ \ \ \ \ \ \ \
G_{ss}
\,\equiv\,
\frac{G_{rs}^2}{G_{rr}},
\]
by successive differentiations and replacement,
taking values at the origin, one convinces oneself
({\em see} also~{\cite{Chen-Merker-2020}}),
that all $F_{i,j}$ with $j \geqslant 2$ express
in terms of the $F_{j',0}$ with $j' \leqslant i+j$ 
and of the $F_{j',1}$ with $j' + 1 \leqslant i + j$.
Here, one obtains $F_{2,2} = 2$, $F_{1,3} = 0$, 
$F_{0,4} = 0$,
and the same for $G$.

\begin{Lemma}
\label{Lm-G40-zero}
One can normalize $G_{4,0} := 0$.
\end{Lemma}

\proof
Indeed, $\eqFG$ gives, with the free parameter $b_2$:
\begin{align}
0
&
\overset{\green{\bf 4,0}}{\,\,=\,\,}
-\,\tfrac{1}{24}\,
a_{1,1}^2\,F_{4,0}
+
\tfrac{1}{24}\,
a_{1,1}^4\,G_{4,0}
+
\tfrac{1}{6}\,
a_{1,1}^3\,a_{2,1}\,
G_{3,1}
+
\tfrac{1}{8}\,
a_{1,1}^2\,a_{2,1}^2
+
\tfrac{1}{4}\,
a_{1,1}^2\,\boxed{b_2},
\notag
\\
0
&
\overset{\green{\bf 3,1}}{\,\,=\,\,}
-\,\tfrac{1}{6}\,
a_{1,1}^2\,
F_{3,1}
+
\tfrac{1}{6}\,
a_{1,1}^3\,
G_{3,1}.
\qedhere
\end{align}
\endproof

Visibly, $G_{3,1} \propto F_{3,1}$ is a relative invariant,
and we have:
\[
u
\,=\,
\tfrac{x^2}{2}
+
\tfrac{x^2y}{2}
+
F_{3,1}\,\tfrac{x^3y}{6}
+
\tfrac{x^2y^2}{2}
+
{\rm O}_{x,y}(5)
\ \ \ \ \ \ \ \
\xrightarrow[{\rule[0pt]{50pt}{0pt}}]{\text{Equivalence}}
\ \ \ \ \ \ \ \
v
\,=\,
\tfrac{r^2}{2}
+
\tfrac{r^2s}{2}
+
G_{3,1}\,\tfrac{r^3s}{6}
+
\tfrac{r^2s^2}{4}
+
{\rm O}_{r,s}(4).
\]

\begin{Lemma}
\label{Lm-Gstab-4}
Stabilization of this order $\leqslant 4$ normalization holds 
if and only if:
\[
\left[
\begin{array}{ccc}
a_{1,1} & \red{\bf 0} & \red{-a_{1,1}a_{2,1}}
\\
a_{2,1} & \red{\bf 1} & b_2
\\
\red{\bf 0} & \red{\bf 0} & \red{a_{1,1}^2}
\end{array}
\right]^{\green{\bf 3}}
\,\,\,\leadsto\,\,\,
\def\arraystretch{1.25}
\left[
\begin{array}{ccc}
a_{1,1} & \red{\bf 0} & \red{-a_{1,1}a_{2,1}}
\\
a_{2,1} & \red{\bf 1} &
\red{-\tfrac{1}{2}a_{2,1}^2-\tfrac{2}{3}a_{1,1}a_{2,1}G_{3,1}}
\\
\red{\bf 0} & \red{\bf 0} & \red{a_{1,1}^2}
\end{array}
\right]^{\green{\bf 4}}.
\]
\end{Lemma}

\proof
After setting $G_{4,0} := 0 =: F_{4,0}$, 
solve $b_2$ in equation 
$\overset{\green{\bf 4,0}}{\,\,=\,\,}$ above.
\endproof

Coming back to order 3, in the infinitesimal counterpart $\eqL$:
\[
\aligned
0
&
\overset{\green{\bf 3,0}}{\,\,=\,\,}
\tfrac{1}{6}\,F_{3,1}\,T_2
+
\tfrac{1}{2}\,A_{2,1}
+
\tfrac{1}{2}\,B_1,
\\
0
&
\overset{\green{\bf 2,1}}{\,\,=\,\,}
\tfrac{1}{2}\,
F_{3,1}\,T_1
+
\tfrac{1}{2}\,
T_2
+
\tfrac{1}{2}\,
A_{2,2},
\endaligned
\]
we normalize:
\[
\aligned
B_1
&
\,:=\,
-\,\tfrac{1}{3}\,
F_{3,1}\,T_2
-
A_{2,1},
\\
A_{2,2}
&
\,=\,
-\,F_{3,1}\,T_1
-
T_2.
\endaligned
\]

Beyond, because $F_{3,1}$ is a relative 
invariant, we must open two branches:
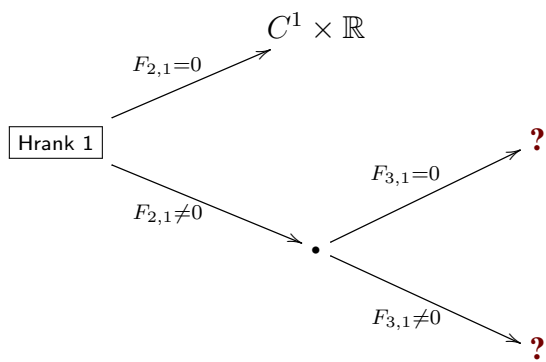
\[
\xymatrix{
&&
C^1\times\R
&&
\\
\ar[urr]^{F_{2,1}=0}
\boxed{\scriptstyle\sf Hrank~1}
\ar[drr]_{F_{2,1}\neq0}
&&
&&
\red{\text{\bf ?}}
\\
&&
\ar[urr]^{F_{3,1}=0}
\centersmallbullet
\ar[drr]_{F_{3,1}\neq 0}
&&
\\
&&
&&
\red{\text{\bf ?}}
}
\]

We study first the branch $F_{3,1} = 0$, and we let terms of
order 5 appear:
\[
u
\,=\,
\tfrac{x^2}{2}
+
\tfrac{x^2y}{2}
+
\tfrac{x^2y^2}{2}
+
F_{5,0}\,\tfrac{x^5}{120}
+
F_{4,1}\,\tfrac{x^4y}{24}
+
\tfrac{x^2y^3}{2}
+
{\rm O}_{x,y}(6).
\]

\begin{Lemma}
In the branch $F_{3,1} = 0$, 
affine homogeneity forces $F_{4,1} = 0$, necessarily.
\end{Lemma}

\proof
Indeed, $\eqL$ gives:
\begin{align}
0
&
\overset{\green{\bf 4,0}}{\,\,=\,\,}
\tfrac{1}{24}\,
F_{5,0}\,T_1
+
\tfrac{1}{24}\,
F_{4,1}\,T_2
+
\tfrac{1}{4}\,
B_2,
\notag
\\
0
&
\overset{\green{\bf 3,1}}{\,\,=\,\,}
\tfrac{1}{6}\,
F_{4,1}\,
T_1.
\qedhere
\end{align}
\endproof

Next, solve:
\[
B_2
\,:=\,
-\,\tfrac{1}{6}\,
F_{5,0}\,T_1
-
0.
\]
Thus with $F_{4,1} := 0 =: G_{4,1}$:
\[
\aligned
u
&
\,=\,
\tfrac{x^2}{2}
+
\tfrac{x^2y}{2}
+
\tfrac{x^2y^2}{2}
+
F_{5,0}\,\tfrac{x^5}{120}
+
0
+
\tfrac{x^2y^3}{2}
+
{\rm O}_{x,y}(6),
\\
\xrightarrow[{\rule[0pt]{50pt}{0pt}}]{\text{Equivalence}}
\ \ \ \ \ \ \ \ \ \ \ \ \ \ \ \ \ \ \ \
v
&
\,=\,
\tfrac{r^2}{2}
+
\tfrac{r^2s}{2}
+
\tfrac{r^2s^2}{2}
+
G_{5,0}\,\tfrac{r^5}{120}
+
0
+
\tfrac{r^2s^3}{2}
+
{\rm O}_{r,s}(6),
\endaligned
\]

With this, $\eqFG$
at order $\mu = 5$ contains only one nonzero equation:
\[
0
\overset{\green{\bf 5,0}}{\,\,=\,\,}
-\,\tfrac{1}{120}\,
F_{5,0}\,a_{1,1}^2
+
\tfrac{1}{120}\,
G_{5,0}\,
a_{1,1}^5.
\]
Therefore, $F_{5,0} \propto G_{5,0}$ is a relative invariant:
it creates a new branching:
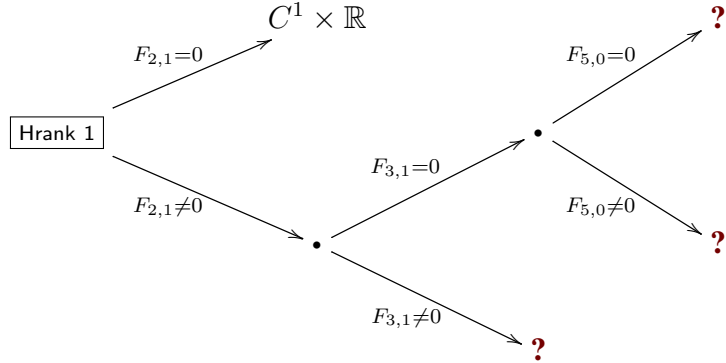
\[
\xymatrix{
&&
C^1\times\R
&&
&&
\red{\text{\bf ?}}
\\
\ar[urr]^{F_{2,1}=0}
\boxed{\scriptstyle\sf Hrank~1}
\ar[drr]_{F_{2,1}\neq0}
&&
&&
\ar[urr]^{F_{5,0}=0}
\centersmallbullet
\ar[drr]_{F_{5,0}\neq0}
&&
\\
&&
\ar[urr]^{F_{3,1}=0}
\centersmallbullet
\ar[drr]_{F_{3,1}\neq 0}
&&
&&
\red{\text{\bf ?}}
\\
&&
&&
\red{\text{\bf ?}}
}
\]

Study first the subbranch $F_{5,0} = 0$:
\[
u
\,=\,
\tfrac{x^2}{2}
+
\tfrac{x^2y}{2}
+
\tfrac{x^2y^2}{2}
+
0
+
0
+
\tfrac{x^2y^3}{2}
+
{\rm O}_{x,y}(6),
\]
with $v = G$ similarly given. 
At order 4, the isotropy group
from Lemma~{\ref{Lm-Gstab-4}}:
\[
\def\arraystretch{1.25}
\left[
\begin{array}{ccc}
a_{1,1} & \red{\bf 0} & \red{-a_{1,1}a_{2,1}}
\\
a_{2,1} & \red{\bf 1} &
\red{-\tfrac{1}{2}a_{2,1}^2}
\\
\red{\bf 0} & \red{\bf 0} & \red{a_{1,1}^2}
\end{array}
\right]^{\green{\bf 4}},
\]
is still 2-dimensional, with parameters $a_{1,1}$, $a_{2,1}$.

\begin{Proposition}
In the branch $F_{2,1} \neq 0$, $F_{3,1} = 0$, $F_{5,0} = 0$,
if the surface $S^2 \subset \C^3$ is affinely homogeneous,
then all $F_{j,k} = 0$ except $F_{2,k} = k!$
for every $k = 1, 2, 3, 4, 5, \dots$.
\end{Proposition}

\proof
Examine $\eqL$ at order 5:
\[
\aligned
0
&
\overset{\green{\bf 5,0}}{\,\,=\,\,}
\tfrac{1}{120}\,
F_{6,0}\,T_1
+
\tfrac{1}{120}\,
F_{5,1}\,T_2,
\\
0
&
\overset{\green{\bf 4,1}}{\,\,=\,\,}
\tfrac{1}{24}\,
F_{5,1}\,T_1
+
\tfrac{1}{24}\,
F_{4,2}\,T_2,
\\
0
&
\overset{\green{\bf 3,2}}{\,\,=\,\,}
\tfrac{1}{12}\,
F_{4,2}\,T_1
+
\tfrac{1}{12}\,
F_{3,3}\,T_2,
\\
0
&
\overset{\green{\bf 2,3}}{\,\,=\,\,}
\tfrac{1}{12}\,
F_{3,3}\,T_1
+
\big(
-2
+
\tfrac{1}{12}\,F_{2,4}
\big)\,
T_2,
\\
0
&
\overset{\green{\bf 1,4}}{\,\,=\,\,}
\big(
\tfrac{1}{24}\,F_{2,4}
-
1
\big)\,
T_1
+
\tfrac{1}{24}\,
F_{1,5}\,
T_2,
\\
0
&
\overset{\green{\bf 0,5}}{\,\,=\,\,}
\tfrac{1}{120}\,
F_{1,5}\,
T_1
+
\tfrac{1}{120}\,
F_{0,6}\,
T_2,
\endaligned
\]
to get:
\[
F_{2,4}
\,=\,
4!
\ \ \ \ \ \ \ \ \ \ \ \ \ \ \ \ \ \ \ \
\text{while}
\ \ \ \ \ \ \ \ \ \ \ \ \ \ \ \ \ \ \ \
0
\,=\,
F_{0,6}
\,=\,
F_{1,5}
\,=\,
F_{3,3}
\,=\,
F_{4,2}
\,=\,
F_{5,1}
\,=\,
F_{6,0}.
\]

Next, $\eqL$ at order 6:
\[
\aligned
0
&
\overset{\green{\bf 6,0}}{\,\,=\,\,}
\tfrac{1}{720}\,F_{7,0}\,T_1
+
\tfrac{1}{720}\,F_{6,1}\,T_2,
\\
0
&
\overset{\green{\bf 5,1}}{\,\,=\,\,}
\tfrac{1}{120}\,F_{6,1}\,T_1
+
\tfrac{1}{120}\,F_{5,2}\,T_2,
\\
0
&
\overset{\green{\bf 4,2}}{\,\,=\,\,}
\tfrac{1}{48}\,F_{5,2}\,T_1
+
\tfrac{1}{48}\,F_{4,3}\,T_2,
\\
0
&
\overset{\green{\bf 3,3}}{\,\,=\,\,}
\tfrac{1}{36}\,F_{4,3}\,T_1
+
\tfrac{1}{36}\,F_{3,4}\,T_2,
\\
0
&
\overset{\green{\bf 2,4}}{\,\,=\,\,}
\tfrac{1}{48}\,F_{3,4}\,T_1
+
\big(
-\tfrac{5}{2}
+
\tfrac{1}{48}\,F_{25}
\big)\,T_2,
\\
0
&
\overset{\green{\bf 1,5}}{\,\,=\,\,}
\big(
\tfrac{1}{120}\,F_{2,5}
-
1
\big)\,T_1
+
\tfrac{1}{120}\,F_{1,6}\,T_1,
\\
0
&
\overset{\green{\bf 0,6}}{\,\,=\,\,}
\tfrac{1}{720}\,F_{1,6}\,T_1
+
\tfrac{1}{720}\,F_{0,7}\,T_2.
\endaligned
\]
solve similarly:
\[
F_{2,5}
\,=\,
5!
\ \ \ \ \ \ \ \ \ \ \ \ \ \ \ 
\text{while}
\ \ \ \ \ \ \ \ \ \ \ \ \ \ \ 
0
\,=\,
F_{0,7}
\,=\,
F_{1,6}
\,=\,
F_{3,4}
\,=\,
F_{4,3}
\,=\,
F_{5,2}
\,=\,
F_{6,1}
\,=\,
F_{7,0}.
\]
An induction on the order $\mu = i + j$ is elementary.
\endproof

Since $\sum_k\, y^k = \frac{1}{1-y}$, we obtain

\begin{Theorem}
In the branch $F_{2,1} \neq 0$,
$F_{3,1} = 0$, $F_{5,0} = 0$, there is a
single affinely homogeneous surface $S^2 \subset \R^3$:
\[
u
\,=\,
\frac{1}{2}\,
\frac{x^2}{1-y},
\]
which has 4-dimensional transitive affine
Lie symmetry algebra generated by:
\[
\aligned
e_1
&
\,:=\,
(1-y)\,\partial_x
+
x\,\partial_u,
\\
e_2
&
\,:=\,
(1-y)\,\partial_y
+
u\,\partial_u,
\\
e_3
&
\,:=\,
x\,\partial_x
+
2\,u\,\partial_u,
\\
e_4
&
\,:=\,
-\,u\,\partial_x
+
x\,\partial_y,
\endaligned
\]
sharing the Lie brackets:
\[
[e_1,e_2]
\,=\,
e_1,
\ \ \ \ \
[e_1,e_3]
\,=\,
e_1,
\ \ \ \ \
[e_1,e_4]
\,=\,
e_2,
\ \ \ \ \
[e_2,e_4]
\,=\,
e_4,
\ \ \ \ \
[e_3,e_4]
\,=\,
e_4.
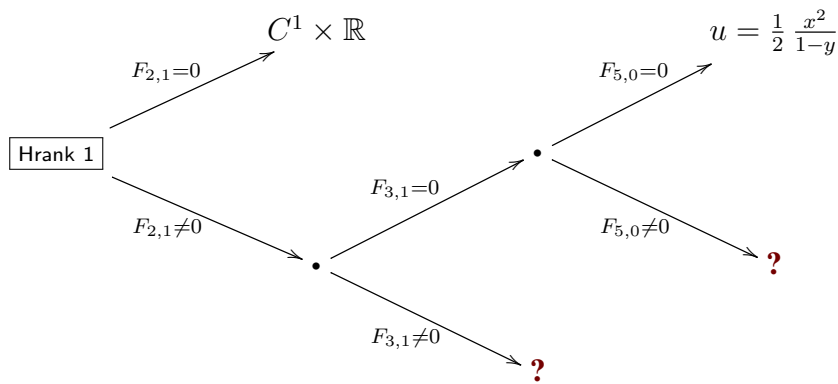
\eqno\qed
\]
\end{Theorem}

\[
\xymatrix{
&&
C^1\times\R
&&
&&
u=\tfrac{1}{2}\,\tfrac{x^2}{1-y}
\\
\ar[urr]^{F_{2,1}=0}
\boxed{\scriptstyle\sf Hrank~1}
\ar[drr]_{F_{2,1}\neq0}
&&
&&
\ar[urr]^{F_{5,0}=0}
\centersmallbullet
\ar[drr]_{F_{5,0}\neq0}
&&
\\
&&
\ar[urr]^{F_{3,1}=0}
\centersmallbullet
\ar[drr]_{F_{3,1}\neq 0}
&&
&&
\red{\text{\bf ?}}
\\
&&
&&
\red{\text{\bf ?}}
}
\]

Next, let us study the subbranch $F_{5,0} \neq 0$. From:
\[
0
\overset{\green{\bf 5,0}}{\,\,=\,\,}
-\,\tfrac{1}{120}\,
F_{5,0}\,a_{1,1}^2
+
\tfrac{1}{120}\,
G_{5,0}\,
a_{1,1}^5,
\]
taking $a_{1,1} := \sqrt[3]{F_{5,0}}$, 
we normalize $G_{5,0} := 1 =: F_{5,0}$.
To stabilize:
\[
u
\,=\,
\tfrac{x^2}{2}
+
\tfrac{x^2y}{2}
+
\tfrac{x^2y^2}{2}
+
\tfrac{x^5}{120}
+
\tfrac{x^2y^3}{2}
+
{\rm O}_{x,y}(6)
\xrightarrow[{\rule[0pt]{50pt}{0pt}}]{\text{Equivalence}}
v
\,=\,
\tfrac{r^2}{2}
+
\tfrac{r^2s}{2}
+
\tfrac{r^2s^2}{2}
+
\tfrac{r^5}{120}
+
\tfrac{r^2s^3}{2}
+
{\rm O}_{r,s}(6),
\]
we need to satisfy:
\[
0
\overset{\green{\bf 5,0}}{\,\,=\,\,}
-\,\tfrac{1}{120}\,
a_{1,1}^2
+
\tfrac{1}{120}\,
a_{1,1}^5,
\]
and we set $a_{1,1} := 1$.

At the infinitesimal level, $\eqL$ for order 5 gives:
\[
\aligned
0
&
\overset{\green{\bf 5,0}}{\,\,=\,\,}
\tfrac{1}{120}\,F_{6,0}\,T_1
+
\big(
-\tfrac{1}{120}
+
\tfrac{1}{120}\,F_{5,1}
\big)\,
T_2
+
\tfrac{1}{40}\,
A_{1,1},
\\
0
&
\overset{\green{\bf 4,1}}{\,\,=\,\,}
\big(
\tfrac{1}{24}\,F_{5,1}
-
\tfrac{1}{6}
\big)\,
T_1
+
\tfrac{1}{24}\,F_{4,2}\,T_2,
\\
0
&
\overset{\green{\bf 3,2}}{\,\,=\,\,}
\tfrac{1}{12}\,
F_{4,2}\,T_1
+
\tfrac{1}{12}\,F_{3,3}\,T_2,
\\
0
&
\overset{\green{\bf 2,3}}{\,\,=\,\,}
\tfrac{1}{12}\,F_{3,3}\,T_1
+
\big(
-2
+
\tfrac{1}{12}\,F_{2,4}
\big)\,
T_2,
\\
0
&
\overset{\green{\bf 1,4}}{\,\,=\,\,}
\big(
\tfrac{1}{24}\,F_{2,4}
-
1
\big)\,T_1
+
\tfrac{1}{24}\,
F_{1,5}\,T_2,
\\
0
&
\overset{\green{\bf 0,5}}{\,\,=\,\,}
\tfrac{1}{120}\,
F_{1,5}\,T_1
+
\tfrac{1}{120}\,
F_{0,6}\,T_2,
\endaligned
\]
whence:
\[
F_{5,1}
\,=\,
4,
\ \ \ \ \
F_{2,4}
\,=\,
4!
\ \ \ \ \ \ \ \ \ \ \ \ \ \ \ \ \ \ \ \
\text{while}
\ \ \ \ \ \ \ \ \ \ \ \ \ \ \ \ \ \ \ \
0
\,=\,
F_{0,6}
\,=\,
F_{1,5}
\,=\,
F_{3,3}
\,=\,
F_{4,2},
\]
and lastly:
\[
A_{1,1}
\,:=\,
-\,
\tfrac{1}{3}\,F_{6,0}\,T_1
-
T_2.
\]

Next, at order 6, putting similarly as always:
\[
G_{5,1}
\,=\,
4,
\ \ \ \ \
G_{2,4}
\,=\,
4!
\ \ \ \ \ \ \ \ \ \ \ \ \ \ \ \ \ \ \ \
\text{while}
\ \ \ \ \ \ \ \ \ \ \ \ \ \ \ \ \ \ \ \
0
\,=\,
G_{0,6}
\,=\,
G_{1,5}
\,=\,
G_{3,3}
\,=\,
G_{4,2},
\]
only one nontrivial equation exists:
\[
0
\overset{\green{\bf 6,0}}{\,\,=\,\,}
-\,\tfrac{1}{720}\,F_{6,0}
+
\tfrac{1}{720}\,G_{6,0}
+
\tfrac{1}{240}\,
\boxed{a_{2,1}}.
\]
Using $a_{2,1}$, we normalize:
\[
G_{6,0}
\,:=\,
0
\,=:\,
F_{6,0},
\]
and then we stabilize:
\[
a_{2,1}
\,:=\,
0.
\]

Since the isotropy matrix is now reduced to the identity:
\[
\left[\!
\begin{array}{ccc}
\red{\bf 1} & \red{\bf 0} & \red{\bf 0}
\\
\red{\bf 0} & \red{\bf 1} & \red{\bf 0}
\\
\red{\bf 0} & \red{\bf 0} & \red{\bf 1}
\end{array}
\!\right],
\]
$\eqFG$ is terminated, and only $\eqL$ must be examined further.

At order 6, $\eqL$ gives:
\[
\aligned
0
&
\overset{\green{\bf 6,0}}{\,\,=\,\,}
\tfrac{1}{720}\,F_{7,0}\,T_1
+
\tfrac{1}{720}\,F_{6,1}\,T_2
+
\tfrac{1}{240}\,A_{2,1},
\\
0
&
\overset{\green{\bf 5,1}}{\,\,=\,\,}
\tfrac{1}{120}\,F_{6,1}\,T_1
+
\big(
-\tfrac{1}{6}
+
\tfrac{1}{120}\,F_{5,2}
\big)\,T_2,
\\
0
&
\overset{\green{\bf 4,2}}{\,\,=\,\,}
\big(
\tfrac{1}{48}\,F_{5,2}
-
\tfrac{5}{12}
\big)\,
T_1
+
\tfrac{1}{48}\,F_{4,3}\,T_2,
\\
0
&
\overset{\green{\bf 3,3}}{\,\,=\,\,}
\tfrac{1}{36}\,F_{4,3}\,T_1
+
\tfrac{1}{36}\,F_{3,4}\,T_2,
\\
0
&
\overset{\green{\bf 2,4}}{\,\,=\,\,}
\tfrac{1}{48}\,F_{3,4}\,T_1
+
\big(
-\tfrac{5}{2}
+
\tfrac{1}{48}\,F_{2,5}
\big)\,T_2,
\\
0
&
\overset{\green{\bf 1,5}}{\,\,=\,\,}
\big(
\tfrac{1}{120}\,F_{2,5}
-
1
\big)\,T_1
+
\tfrac{1}{120}\,F_{1,6}\,T_2,
\\
0
&
\overset{\green{\bf 0,6}}{\,\,=\,\,}
\tfrac{1}{720}\,F_{1,6}\,T_1
+
\tfrac{1}{720}\,F_{0,7}\,T_2,
\endaligned
\]
whence:
\[
F_{0,7}
\,:=\,
0,
\ \ \ \ \
F_{1,6}
\,:=\,
0,
\ \ \ \ \
F_{2,5}
\,:=\,
120,
\ \ \ \ \
F_{3,4}
\,:=\,
0,
\ \ \ \ \
F_{4,3}
\,:=\,
0,
\ \ \ \ \
F_{5,2}
\,:=\,
20,
\ \ \ \ \
F_{6,1}
\,:=\,
0,
\]
and lastly:
\[
A_{2,1}
\,:=\,
-\,\tfrac{1}{3}\,
F_{7,0}\,
T_1.
\]

At order 7, $\eqL$ gives:
\[
\aligned
0
&
\overset{\green{\bf 7,0}}{\,\,=\,\,}
\big(
\tfrac{1}{5040}\,F_{8,0}
-
\tfrac{1}{288}
\big)\,T_1
+
\big(
-\tfrac{1}{840}\,F_{7,0}
+
\tfrac{1}{5040}\,F_{7,1}
\big)\,T_2,
\\
0
&
\overset{\green{\bf 6,1}}{\,\,=\,\,}
\big(
\tfrac{1}{720}\,F_{7,1}
-
\tfrac{1}{120}\,F_{7,0}
\big)\,T_1
+
\tfrac{1}{720}\,F_{6,2}\,T_2,
\\
0
&
\overset{\green{\bf 5,2}}{\,\,=\,\,}
\tfrac{1}{240}\,F_{6,2}\,T_1
+
\big(
-\tfrac{1}{2}
+
\tfrac{1}{240}\,F_{5,3}
\big)\,T_2,
\\
0
&
\overset{\green{\bf 4,3}}{\,\,=\,\,}
\big(
\tfrac{1}{144}\,F_{5,3}
-
\tfrac{5}{6}
\big)\,T_1
+
\tfrac{1}{144}\,F_{4,4}\,T_2,
\\
0
&
\overset{\green{\bf 3,4}}{\,\,=\,\,}
\tfrac{1}{144}\,F_{4,4}\,T_1
+
\tfrac{1}{144}\,F_{3,5}\,T_2,
\\
0
&
\overset{\green{\bf 2,5}}{\,\,=\,\,}
\tfrac{1}{240}\,F_{3,5}\,T_1
+
\big(
-3+\tfrac{1}{240}\,F_{2,6}
\big)\,T_2,
\\
0
&
\overset{\green{\bf 1,6}}{\,\,=\,\,}
\big(
\tfrac{1}{720}\,F_{2,6}
-
1
\big)\,
T_1
+
\tfrac{1}{720}\,F_{1,7}\,T_2,
\\
0
&
\overset{\green{\bf 0,7}}{\,\,=\,\,}
\tfrac{1}{5040}\,F_{1,7}\,T_1
+
\tfrac{1}{5040}\,F_{0,8}\,T_2,
\endaligned
\]
which is solved as:
\[
\aligned
F_{0,8}
\,:=\,
0,
\ \ \ \ \
F_{1,7}
\,:=\,
0,
\ \ \ \ \
F_{2,6}
\,:=\,
720,
\ \ \ \ \
F_{3,5}
&
\,:=\,
0,
\ \ \ \ \
F_{4,4}
\,:=\,
0,
\ \ \ \ \
F_{5,3}
\,:=\,
120,
\\
F_{6,2}
&
\,:=\,
0,
\ \ \ \ \
F_{7,1}
\,:=\,
6\,F_{7,0},
\ \ \ \ \
F_{8,0}
\,:=\,
\tfrac{35}{2},
\endaligned
\]
with an invariant:
\[
F_{7,0}
\,=:\,
\theta
\,\in\,
\R,
\]
which may take any real value.

At order 8, the resolution of the (unwritten)
equations of $\eqL$ is:
\[
\aligned
F_{0,9}
\,:=\,
0,
\ \ \ \ \
F_{1,8}
\,:=\,
0,
\ \ \ \ \
F_{2,7}
\,:=\,
5040,
\ \ \ \ \
F_{3,6}
\,:=\,
0,
\ \ \ \ \
F_{4,5}
&
\,:=\,
0,
\ \ \ \ \
F_{5,4}
\,:=\,
840,
\ \ \ \ \
F_{6,3}
\,:=\,
0,
\\
F_{7,2}
&
\,:=\,
42\,\theta,
\ \ \ \ \
F_{8,1}
\,:=\,
\tfrac{245}{2},
\ \ \ \ \
F_{9,0}
\,:=\,
4\,\theta^2.
\endaligned
\]

One therefore finds a $1$-parameter family of affinely
inequivalent homogeneous models $\big( S_\theta^2 \big)_{\theta 
\in \R}$.

\begin{Proposition}
In the branch $F_{2,1} \neq 0$, $F_{3,1} = 0$, $F_{5,0} \neq 0$,
there is a 1-parameter family of inequivalent 
affinely homogeneous surfaces $S_\theta^2 \subset \R^3$:
\[
\aligned
u
&
\,=\,
\tfrac{1}{2}\,x^2
+
\tfrac{1}{2}\,x^2y
+
\tfrac{1}{2}\,x^2y^2
\\
&
\ \ \ \ \
+
\tfrac{1}{120}\,
x^5
+
\tfrac{1}{2}\,
x^2y^3
\\
&
\ \ \ \ \
+
\tfrac{1}{30}\,
x^5y
+
\tfrac{1}{2}\,
x^2y^4
\\
&
\ \ \ \ \
+
\tfrac{1}{5040}\,
\theta\,
x^7
+
\tfrac{1}{12}\,
x^5y^2
+
\tfrac{1}{2}\,
x^2y^5
\\
&
\ \ \ \ \
+
\tfrac{1}{2304}\,x^8
+
\tfrac{1}{840}\,\theta\,x^7y
+
\tfrac{1}{6}\,x^5y^3
+
\tfrac{1}{2}\,x^2y^6
\\
&
\ \ \ \ \
+
\tfrac{1}{90720}\,\theta^2\,x^9
+
\tfrac{7}{2304}\,x^8y
+
\tfrac{1}{240}\,\theta\,x^7y^2
+
\tfrac{7}{24}\,x^5y^4
+
\tfrac{1}{2}\,x^2y^7
+
{\rm O}_{x,y}(10).
\endaligned
\]
with $2$-dimensional 
(simply transitive) commutative affine Lie symmetry algebra:
\reqnomode\usetagform{EngelLie}
\begin{align}
e_1
&
\,:=\,
\big(
1-y+\tfrac{1}{3}\,\theta\,u
\big)\,
\partial_x
+
\big(
-\,\tfrac{1}{3}\,\theta\,x
-\tfrac{1}{6}\,u
\big)\,
\partial_y
+
x\,\partial_u,
\notag
\\
e_2
&
\,:=\,
-\,x\,\partial_x
+
(1-y)\,\partial_y
-
u\,\partial_u,
\ \ \ \ \ \ \ \ \ \ \ \ \ \ \ \ \ \ \ \ \ \ \ \ \ \ \ \ \ \ \ \ \ \ \ 
\ \ \ \ \ \ \ \ \ \
[e_1,e_2]
\,=\,
0.
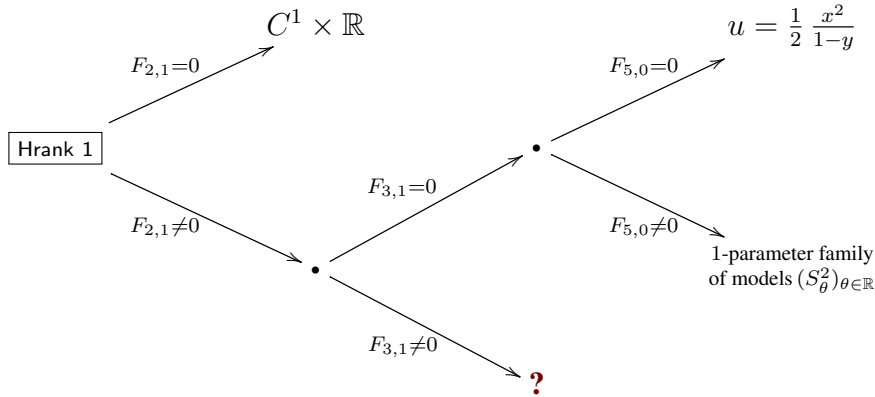
\tag{\qed}
\end{align}
\end{Proposition}

\[
\xymatrix{
&&
C^1\times\R
&&
&&
u=\tfrac{1}{2}\,\tfrac{x^2}{1-y}
\\
\ar[urr]^{F_{2,1}=0}
\boxed{\scriptstyle\sf Hrank~1}
\ar[drr]_{F_{2,1}\neq0}
&&
&&
\ar[urr]^{F_{5,0}=0}
\centersmallbullet
\ar[drr]_{F_{5,0}\neq0}
&&
\\
&&
\ar[urr]^{F_{3,1}=0}
\centersmallbullet
\ar[drr]_{F_{3,1}\neq 0}
&&
&&
{\substack{
1\text{-parameter family}
\\
\text{of models}\,(S_\theta^2)_{\theta\in\R}}}
\\
&&
&&
\red{\text{\bf ?}}
}
\]

It remains to explore the subbranch $F_{3,1} \neq 0$,
within the branch $F_{2,1} = 1$.
From the proof of Lemma~{\ref{Lm-G40-zero}}:
\[
0
\overset{\green{\bf 3,1}}{\,\,=\,\,}
-\,\tfrac{1}{6}\,a_{1,1}^2\,F_{3,1}
+
\tfrac{1}{6}\,a_{1,1}^3\,G_{3,1},
\]
it is clear that we can normalize $G_{3,1} := 1 =: F_{3,1}$.
Thus:
\[
u
\,=\,
\tfrac{x^2}{2}
+
\tfrac{x^2y}{2}
+
\tfrac{x^3y}{6}
+
\tfrac{x^2y^2}{2}
+
{\rm O}_{x,y}(5)
\ \ \ \ \ \ \ \
\xrightarrow[{\rule[0pt]{50pt}{0pt}}]{\text{Equivalence}}
\ \ \ \ \ \ \ \
v
\,=\,
\tfrac{r^2}{2}
+
\tfrac{r^2s}{2}
+
\tfrac{r^3s}{6}
+
\tfrac{r^2s^2}{2}
+
{\rm O}_{r,s}(5).
\]

Since $\overset{\green{\bf 3,1}}{\,=\,}$ becomes
$0 = -\,\frac{1}{6}\, a_{1,1}^2 + \frac{1}{6}\, a_{1,1}^3$,
we set $a_{1,1} := 1$, hence 
the stability group at orders $\leqslant 4$ is 1-dimensional:
\[
\def\arraystretch{1.25}
\left[
\begin{array}{ccc}
\red{\bf 1} & \red{\bf 0} & \red{-a_{2,1}}
\\
a_{2,1} & \red{\bf 1} &
\red{-\tfrac{1}{2}a_{2,1}^2-\tfrac{2}{3}a_{2,1}}
\\
\red{\bf 0} & \red{\bf 0} & \red{a_{1,1}^2}
\end{array}
\right]^{\green{\bf 4}}.
\]

Next, $\eqL$ at order $4$ gives:
\[
\aligned
0
&
\overset{\green{\bf 4,0}}{\,\,=\,\,}
\tfrac{1}{24}\,F_{5,0}\,T_1
+
\tfrac{1}{24}\,F_{4,1}\,T_2
+
\tfrac{1}{6}\,A_{2,1}
+
\tfrac{1}{4}\,B_2,
\\
0
&
\overset{\green{\bf 3,1}}{\,\,=\,\,}
\big(
\tfrac{1}{6}\,F_{4,1}
-
\tfrac{1}{6}
\big)\,T_1
+
\big(
-\tfrac{2}{3}
+
\tfrac{1}{6}\,F_{3,2}
\big)\,T_2,
\\
0
&
\overset{\green{\bf 2,2}}{\,\,=\,\,}
\big(
\tfrac{1}{4}\,F_{3,2}
-
\tfrac{3}{2}
\big)\,T_1
+
\big(
-\tfrac{3}{2}
+
\tfrac{1}{4}\,F_{2,3}
\big)\,T_2,
\\
0
&
\overset{\green{\bf 1,3}}{\,\,=\,\,}
\big(
\tfrac{1}{6}\,F_{2,3}
-
1
\big)\,T_1
+
\tfrac{1}{6}\,F_{1,4}\,T_2,
\\
0
&
\overset{\green{\bf 0,4}}{\,\,=\,\,}
\tfrac{1}{24}\,F_{1,4}\,T_1
+
\tfrac{1}{24}\,F_{0,5}\,T_2,
\endaligned
\]
whence:
\[
0
\,=\,
F_{0,5}
\,=\,
F_{1,4},
\ \ \ \ \ \ \ \ \ \ \ \ \ \ \ \ \ \ \ \ \ \ \ \ \ \
F_{2,3}
\,=\,
6,
\ \ \ \ \ \ \
F_{3,2}
\,=\,
6,
\]
and:
\[
\aligned
A_{1,1}
&
\,:=\,
\big(
-F_{4,1}+1
\big)\,T_1
-
2\,T_2,
\\
B_2
&
\,:=\,
-\,\tfrac{1}{6}\,F_{5,0}\,T_1
-
\tfrac{1}{6}\,F_{4,1}\,T_2
-
\tfrac{2}{3}\,A_{2,1}.
\endaligned
\]

Putting in $\eqFG$ at order 5:
\[
0
\,=\,
G_{0,5}
\,=\,
G_{1,4},
\ \ \ \ \ \ \ \ \ \ \ \ \ \ \ \ \ \ \ \ \ \ \ \ \ \
G_{2,3}
\,=\,
6,
\ \ \ \ \ \ \
G_{3,2}
\,=\,
6,
\]
we get:
\[
\aligned
0
&
\overset{\green{\bf 5,0}}{\,\,=\,\,}
-\,\tfrac{1}{120}\,F_{5,0}
+
\tfrac{1}{120}\,G_{5,0}
+
\tfrac{1}{24}\,G_{4,1}\,a_{2,1}
-
\tfrac{1}{18}\,a_{2,1}
+
\tfrac{1}{24}\,a_{2,1}^2,
\\
0
&
\overset{\green{\bf 4,1}}{\,\,=\,\,}
-\,\tfrac{1}{24}\,F_{4,1}
+
\tfrac{1}{24}\,G_{4,1}
+
\tfrac{1}{12}\,\boxed{a_{2,1}}.
\endaligned
\]

Using the last remaining group parameter $a_{2,1}$, we
normalize $G_{4,1} := 0 =: F_{4,1}$, whence $a_{2,1} := 0$
so that the group reduction descends to identity:
\[
\def\arraystretch{1.25}
\left[
\begin{array}{ccc}
\red{\bf 1} & \red{\bf 0} & \red{-a_{2,1}}
\\
a_{2,1} & \red{\bf 1} &
\red{-\tfrac{1}{2}a_{2,1}^2-\tfrac{2}{3}a_{2,1}}
\\
\red{\bf 0} & \red{\bf 0} & \red{a_{1,1}^2}
\end{array}
\right]^{\green{\bf 4}}
\,\,\,\leadsto\,\,\,
\left[\!
\begin{array}{ccc}
\red{\bf 1} & \red{\bf 0} & \red{\bf 0}
\\
\red{\bf 0} & \red{\bf 1} & \red{\bf 0}
\\
\red{\bf 0} & \red{\bf 0} & \red{\bf 1}
\end{array}
\!\right]^{\green{\bf 5}},
\]
to stabilize the normal forms at orders $\leqslant 5$:
\[
\aligned
u
&
\,=\,
\tfrac{x^2}{2}
+
\tfrac{x^2y}{2}
+
\tfrac{x^3y}{6}
+
\tfrac{x^2y^2}{2}
+
\tfrac{1}{120}\,F_{5,0}\,x^5
+
\tfrac{1}{2}\,x^3y^2
+
\tfrac{1}{2}\,x^2y^3
+
{\rm O}_{x,y}(6)
\\
\ \ \ \ \ \ \ \
\xrightarrow[{\rule[0pt]{50pt}{0pt}}]{\text{Equivalence}}
\ \ \ \ \ \ \ \
v
&
\,=\,
\tfrac{r^2}{2}
+
\tfrac{r^2s}{2}
+
\tfrac{r^3s}{6}
+
\tfrac{r^2s^2}{2}
+
\tfrac{1}{120}\,G_{5,0}\,r^5
+
\tfrac{1}{2}\,r^3s^2
+
\tfrac{1}{2}\,s^2r^3
+
{\rm O}_{r,s}(6).
\endaligned
\]

Thus, $\eqFG$ is terminated, and only simply transitive 
homogeneous models can be found in this branch.
The only remaining equation:
\[
0
\overset{\green{\bf 5,0}}{\,\,=\,\,}
-\,\tfrac{1}{120}\,F_{5,0}
+
\tfrac{1}{120}\,G_{5,0},
\]
seems to say that $F_{5,0}$ is an absolute invariant
which may take any value $\eta \in \R$, 
but we will see at higher orders
that only 1 specific numeric 
value is possible for $F_{5,0}$.

Next, go to $\eqL$ at order 5:
\[
\aligned
0
&
\overset{\green{\bf 5,0}}{\,\,=\,\,}
\big(
\tfrac{1}{90}\,F_{5,0}
+
\tfrac{1}{120}\,F_{6,0}
\big)\,T_1
+
\big(
-\tfrac{7}{120}\,F_{5,0}
+
\tfrac{1}{120}\,F_{5,1}
\big)\,T_2
-
\tfrac{1}{18}\,
A_{2,1},
\\
0
&
\overset{\green{\bf 4,1}}{\,\,=\,\,}
\big(
\tfrac{1}{24}\,F_{5,1}
-
\tfrac{1}{6}\,F_{5,0}
\big)\,T_1
+
\big(
-\tfrac{5}{36}
+
\tfrac{1}{24}\,F_{4,2}
\big)\,T_2
+
\tfrac{1}{12}\,A_{2,1},
\\
0
&
\overset{\green{\bf 3,2}}{\,\,=\,\,}
\big(
-\tfrac{1}{2}
+
\tfrac{1}{12}\,F_{4,2}
\big)\,T_1
+
\big(
-3
+
\tfrac{1}{12}\,F_{3,3}
\big)\,T_2,
\\
0
&
\overset{\green{\bf 2,3}}{\,\,=\,\,}
\big(
-3
+
\tfrac{1}{12}\,F_{3,3}
\big)\,T_1
+
\big(
-2
+
\tfrac{1}{12}\,F_{2,4}
\big)\,T_2,
\\
0
&
\overset{\green{\bf 1,4}}{\,\,=\,\,}
\big(
\tfrac{1}{24}\,F_{2,4}
-
1
\big)\,T_1
+
\tfrac{1}{24}\,F_{1,5}\,T_2,
\\
0
&
\overset{\green{\bf 0,5}}{\,\,=\,\,}
\tfrac{1}{120}\,F_{1,5}\,T_1
+
\tfrac{1}{120}\,F_{0,6}\,T_2.
\endaligned
\]

Firstly, solve the last four equations:
\[
F_{0,6}
\,:=\,
0,
\ \ \ \ \ \ \
F_{1,5}
\,:=\,
0,
\ \ \ \ \ \ \
F_{2,4}
\,:=\,
24,
\ \ \ \ \ \ \
F_{3,3}
\,:=\,
36,
\ \ \ \ \ \ \
F_{4,2}
\,:=\,
6,
\]
and secondly, solve $\overset{\green{\bf 4,1}}{\,=\,}$:
\[
A_{2,1}
\,:=\,
\big(
-\tfrac{1}{2}\,F_{5,1}
+
2\,F_{5,0}
\big)\,T_1
-
\tfrac{4}{3}\,T_2.
\]

There remains one equation:
\[
0
\overset{\green{\bf 5,0}}{\,\,=\,\,}
\big(
-\,\tfrac{1}{10}\,F_{5,0}
+
\tfrac{1}{120}\,F_{6,0}
+
\tfrac{1}{36}\,F_{5,1}
\big)\,T_1
+
\big(
-\tfrac{7}{120}\,F_{5,0}
+
\tfrac{1}{120}\,F_{5,1}
+
\tfrac{2}{27}
\big)\,T_2.
\]
Since $\{T_1, T_2\}$ must be free, we deduce:
\[
\aligned
0
&
\,=\,
-\,\tfrac{1}{10}\,F_{5,0}
+
\tfrac{1}{120}\,F_{6,0}
+
\tfrac{1}{36}\,F_{5,1},
\\
0
&
\,=\,
-\tfrac{7}{120}\,F_{5,0}
+
\tfrac{1}{120}\,F_{5,1}
+
\tfrac{2}{27},
\endaligned
\]
which we solve by assigning specific values to two
Taylor coefficients of order 6:
\[
\aligned
F_{5,1}
&
\,:=\,
7\,F_{5,0}
-
\tfrac{80}{9},
\\
F_{6,0}
&
\,:=\,
-\,\tfrac{34}{3}\,F_{5,0}
+
\tfrac{800}{27}.
\endaligned
\]
Up to this point, $F_{5,0}$ is still free, and could
be any real number $\eta \in \R$.

Next, from $\eqL$ at order 6:
\[
F_{0,7}
\,:=\,
0,
\ \ \ \ \
F_{1,6}
\,:=\,
0,
\ \ \ \ \
F_{2,5}
\,:=\,
120,
\ \ \ \ \
F_{3,4}
\,:=\,
240,
\ \ \ \ \
F_{4,3}
\,:=\,
90,
\ \ \ \ \
F_{5,2}
\,:=\,
50\,F_{5,0}
-
\tfrac{800}{9},
\]
and it remains:
\[
\aligned
0
&
\overset{\green{\bf 6,0}}{\,\,=\,\,}
\big(
\tfrac{1}{720}\,F_{7,0}
-
\tfrac{40}{243}
-
\tfrac{7}{160}\,F_{5,0}^2
+
\tfrac{8}{45}\,F_{5,0}
\big)\,
T_1
+
\big(
\tfrac{1}{720}\,F_{6,1}
-
\tfrac{22}{81}
+
\tfrac{67}{720}\,F_{5,0}
\big)\,
T_2,
\\
0
&
\overset{\green{\bf 5,1}}{\,\,=\,\,}
\big(
\tfrac{47}{120}\,F_{5,0}
-
\tfrac{34}{27}
+
\tfrac{1}{120}\,F_{6,1}
\big)\,
T_1
+
\big(
-\tfrac{1}{20}\,F_{5,0}
+
\tfrac{1}{9}
\big)\,
T_2.
\endaligned
\]
Surprisingly:
\[
F_{5,0}
\,:=\,
\frac{20}{9}.
\]
Then:
\[
F_{6,1}
\,:=\,
\tfrac{140}{3},
\ \ \ \ \ \ \ \ \ \ \ \ \ \ \ \ \ \ \ \
F_{7,0}
\,:=\,
-\,\tfrac{280}{27}.
\]

\begin{Assertion}
All higher order $F_{j,k}$ with $j + k \geqslant 8$ are
uniquely determined as specific constants.\qed
\end{Assertion}

The infinitesimal symmetries are:
\[
\aligned
L
&
\,=\,
\Big(
\big[
x-y-\tfrac{10}{9}\,u+1
\big]\,T_1
+
\big[
u-2\,x
\big]\,T_2
\Big)\,\tfrac{\partial}{\partial x}
\\
&
\ \ \ \ \
+
\big(
\big[
\tfrac{10}{9}\,x-y-\tfrac{10}{9}\,u
\big]\,T_1
+
\big[
-\tfrac{4}{3}\,x-y+\tfrac{8}{9}\,u+1
\big]\,T_2
\big)\,\tfrac{\partial}{\partial y}
\\
&
\ \ \ \ \
+
\big(
\big[
x+2\,u
\big]\,T_1
+
\big[
-3\,u
\big]\,T_2
\big)\,\tfrac{\partial}{\partial u}.
\endaligned
\]

\begin{Proposition}
In the branch $F_{2,1} \neq 0$, $F_{3,1} \neq 0$, there is a single
affinely homogeneous model:
\[
\aligned
u
&
\,=\,
\tfrac{x^2}{2}
\\
&
\ \ \ \ \ 
+
\tfrac{x^2y}{2}
\\
&
\ \ \ \ \ 
+
\tfrac{x^3y}{6}
+
\tfrac{x^2y^2}{2}
\\
&
\ \ \ \ \ 
+
\tfrac{1}{54}\,x^5
+
\tfrac{1}{2}\,x^3y^2
+
\tfrac{1}{2}\,x^2y^3
\\
&
\ \ \ \ \ 
+
\tfrac{1}{162}\,x^6
+
\tfrac{1}{18}\,x^5y
+
\tfrac{1}{8}\,x^4y^2
+
x^3y^3
+
\tfrac{1}{2}\,x^2y^4
\\
&
\ \ \ \ \ 
-\,\tfrac{1}{486}\,x^7
+
\tfrac{7}{108}\,x^6y
+
\tfrac{5}{54}\,x^5y^2
+
\tfrac{5}{8}\,x^4y^3
+
\tfrac{5}{3}\,x^3y^4
+
\tfrac{1}{2}\,x^2y^5
\\
&
\ \ \ \ \ 
+
\tfrac{5}{5832}\,x^8
+
\tfrac{1}{162}\,x^7y
+
\tfrac{1}{4}\,x^6y^2
+
\tfrac{47}{216}\,x^5y^3
+
\tfrac{15}{8}\,x^4y^4
+
\tfrac{5}{2}\,x^3y^5
+
\tfrac{1}{2}\,x^2y^6
+
{\rm O}_{x,y}(9),
\endaligned
\]
with 2-dimensional (simply-transitive) commutative affine Lie 
symmetry algebra:
generated by:
\[
\aligned
e_1
&
\,:=\,
\big(
x-y-\tfrac{10}{9}\,u+1
\big)\,
\partial_x
+
\big(
\tfrac{10}{9}\,x-y-\tfrac{10}{9}\,u
\big)\,
\partial_y
+
\big(
x+2\,u
\big)\,
\partial_u,
\\
e_2
&
\,:=\,
\big(
u-2\,x
\big)\,
\partial_x
+
\big(
\tfrac{4}{3}\,x-y+\tfrac{8}{9}\,u+1
\big)\,
\partial_y
-
3\,u\,
\partial_u,
\endaligned
\]
having Lie bracket:
\[
[e_1,e_2]
\,=\,
-\,e_1
-
\tfrac{1}{3}\,e_2.
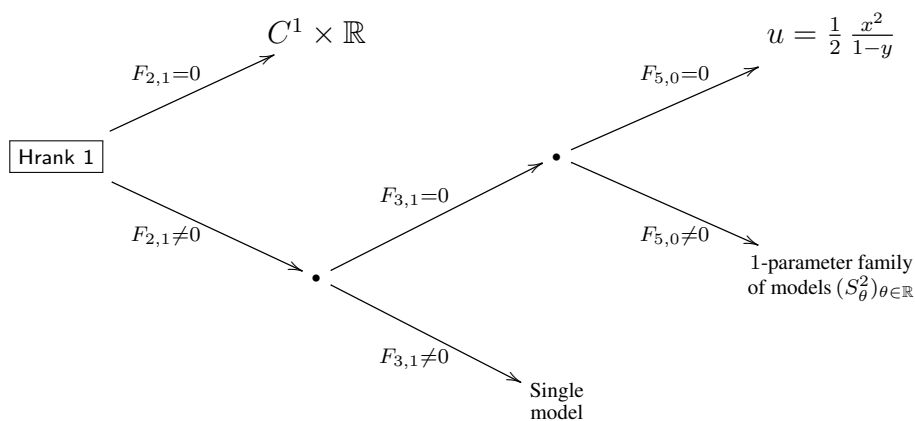
\eqno\qed
\]
\end{Proposition}

\[
\xymatrix{
&&
C^1\times\R
&&
&&
u=\tfrac{1}{2}\,\tfrac{x^2}{1-y}
\\
\ar[urr]^{F_{2,1}=0}
\boxed{\scriptstyle\sf Hrank~1}
\ar[drr]_{F_{2,1}\neq0}
&&
&&
\ar[urr]^{F_{5,0}=0}
\centersmallbullet
\ar[drr]_{F_{5,0}\neq0}
&&
\\
&&
\ar[urr]^{F_{3,1}=0}
\centersmallbullet
\ar[drr]_{F_{3,1}\neq 0}
&&
&&
{\substack{
1\text{-parameter family}
\\
\text{of models}\,(S_\theta^2)_{\theta\in\R}}}
\\
&&
&&
{\substack{
\text{Single}
\\
\text{model}}}
}
\]

\SectionHead{Threefolds $H^3 \subset \R^4$}
{threefolds-H3-R4}

In $\R^4$, consider an affine-linear map
$(x,y,z,u) \longmapsto (r,s,t,v)$ fixing the origin:
\[
\aligned
r
&
\,:=\,
a_{1,1}\,x+a_{1,2}\,y+a_{1,3}\,z+b_1\,u,
\\
s
&
\,:=\,
a_{2,1}\,x+a_{2,2}\,y+a_{2,3}\,z+b_2\,u,
\\
t
&
\,:=\,
a_{3,1}\,x+a_{3,2}\,y+a_{3,3}\,z+b_3\,u,
\\
v
&
\,:=\,
c_1\,x+c_2\,y+c_3\,z+d\,u,
\endaligned
\ \ \ \ \ \ \ \ \ \ \ \ \ \ \ \ \ \ \ \
\text{with}
\ \ \ \ \ \ \ \ \ \ \ \ \ \ \ \ \ \ \ \
0
\,\neq\,
\left\vert\!
\begin{array}{cccc}
a_{1,1} & a_{1,2} & a_{1,3} & b_1
\\
a_{2,1} & b_{2,2} & a_{2,3} & b_2
\\
a_{3,1} & b_{3,2} & a_{3,3} & b_3
\\
c_1 & c_2 & c_3 & d
\end{array}
\!\right\vert.
\]
Also, consider two graphed analytic hypersurfaces:
\[
u
\,=\,
F(x,y,z)
\ \ \ \ \ \ \ \ 
{\scriptstyle{(F(0,0,0)\,=\,0)}}
\ \ \ \ \ \ \ \ \ \ \ \
\text{and}
\ \ \ \ \ \ \ \ \ \ \ \ 
v
\,=\,
G(r,s,t)
\ \ \ \ \ \ \ \ 
{\scriptstyle{(0\,=\,G(0,0,0))}},
\]
and assume that the above map is an affine equivalence 
$\{ u = F\} \longrightarrow \{ v = G\}$.

As in~{\cite{Merker-2022}}, the main hypothesis of
constant Hessian rank 1, after elementary preliminary
transformations:
\[
u
\,=\,
\tfrac{x^2}{2}
+
{\rm O}_{x,y,z}(3),
\]
reads as:
\[
\aligned
1
\,\equiv\,
\rank\,
\left[\!
\begin{array}{ccc}
F_{xx} & F_{xy} & F_{xz}
\\
F_{yx} & F_{yy} & F_{yz}
\\
F_{zx} & F_{zy} & F_{zz}
\end{array}
\!\right],
\endaligned
\]
which is then equivalent to:
\[
0
\,\equiv\,
\left\vert\!
\begin{array}{cc}
F_{xx} & F_{xy} 
\\
F_{yx} & F_{yy}
\end{array}
\!\right\vert
\,\equiv\,
\left\vert\!
\begin{array}{cc}
F_{xx} & F_{xz} 
\\
F_{yx} & F_{yz}
\end{array}
\!\right\vert
\,\equiv\,
\left\vert\!
\begin{array}{cc}
F_{xx} & F_{xy} 
\\
F_{zx} & F_{zy}
\end{array}
\!\right\vert
\,\equiv\,
\left\vert\!
\begin{array}{cc}
F_{xx} & F_{xz} 
\\
F_{zx} & F_{zz}
\end{array}
\!\right\vert.
\]
By affine invariancy of the Hessian matrix rank, 
the same holds about $v = \tfrac{r^2}{2} + {\rm O}_{r,s,t}(3)$.

The fundamental equation which 
holds identically in $\R\{x,y,z\}$:
\[
\aligned
0
\,\equiv\,
\eqFG(x,y,z),
\endaligned
\]
writes:
\[
\aligned
\eqFG
&
\,:=\,
-\,c_1\,x-c_2\,y-c_3\,z-d\,F(x,y,z)
\\
&
\ \ \ \ \
+
G
\Big(
a_{1,1}x+a_{1,2}y+a_{1,3}z+b_1F(x,y,z),\,\,
a_{2,1}x+a_{2,2}y+a_{2,3}z+b_2F(x,y,z),\,\,
\\
&
\ \ \ \ \ \ \ \ \ \ \ \ \ \ \ \ \ \ \ \ \ \ \ \ \ \ \ \ \ \ \ \ \ \ \
\ \ \ \ \ \ \ \ \ \ \ \ \ \ \ \ \ \ \ \ \ \ \ \ \ \ \ \ \ \ \ \ \ \ \
\ \ \ \ \ \
a_{3,1}x+a_{3,2}y+a_{3,3}z+b_3F(x,y,z)
\Big).
\endaligned
\]

Also, an affine vector field:
\[
\aligned
L
&
\,=\,
\ \ 
\big(
T_1+A_{1,1}\,x+A_{1,2}\,y+A_{1,3}\,z+B_1\,u
\big)\,\frac{\partial}{\partial x}
\\
&
\ \ \ \ \
+
\big(
T_2+A_{2,1}\,x+A_{2,2}\,y+A_{2,3}\,z+B_2\,u
\big)\,\frac{\partial}{\partial y}
\\
&
\ \ \ \ \
+
\big(
T_3+A_{3,1}\,x+A_{3,2}\,y+A_{3,3}\,z+B_3\,u
\big)\,\frac{\partial}{\partial z}
\\
&
\ \ \ \ \
+
\big(
T_0+C_1\,x+C_2\,y+C_3\,z+D\,u
\big)\,\frac{\partial}{\partial u},
\endaligned
\]
is tangent to $\{u = F(x,y,z)\}$ if and only if:
\[
\aligned
0
&
\,\equiv\,
\eqL(x,y,z)
\\
&
\,=:\,
L
\big(
-\,u+F(x,y,z)
\big)
\Big\vert_{u=F(x,y,z)},
\endaligned
\]
identically as power series in $\R\{x,y,z\}$.

According to Theorems~1.4, 13.1, 1.5, 25.2
in~{\cite{Merker-2022}}, 
if $H^3 \subset \R^4$ is not affinely equivalent to a product
with $\R^1$ or with $\R^2$, 
its graphing function $F(x,y,z)$ can be 
{\sl pre-normalized}\,\,---\,\,that is, 
{\em normalized before creating any branching}\,\,---\,\,up 
to order $3 + 5 = 8$ included 
and modulo ${\rm O}_{y,z}(3)$ as:
\[
\footnotesize
\aligned
u
&
\,=\,
\ \
\tfrac{x^2}{2}
\\
&
\ \ \ \ \
+
\tfrac{x^2y}{2}
\\
&
\ \ \ \ \
+
\tfrac{x^3z}{6}
+
\tfrac{x^2y^2}{2}
\\
&
\ \ \ \ \
+
F_{4,1,0}\,
\tfrac{x^4y}{24}
+
\tfrac{x^3yz}{2}
\\
&
\ \ \ \ \
+
F_{6,0,0}\,
\tfrac{x^6}{720}
+
F_{5,1,0}\,
\tfrac{x^5y}{120}
+
F_{4,1,0}\,
\tfrac{x^4y^2}{6}
+
\tfrac{x^4z^2}{8}
\\
&
\ \ \ \ \
+
F_{7,0,0}\,
\tfrac{x^7}{5040}
+
F_{6,1,0}\,
\tfrac{x^6y}{720}
+
F_{6,0,1}\,
\tfrac{x^6z}{720}
+
F_{5,1,0}\,
\tfrac{x^5y^2}{24}
+
F_{4,1,0}\,
\tfrac{x^5yz}{12}
\\
&
\ \ \ \ \
+
F_{8,0,0}\,
\tfrac{x^8}{40320}
+
F_{7,1,0}\,
\tfrac{x^7y}{5040}
+
F_{7,0,1}\,
\tfrac{x^7z}{5040}
+
\big(
\tfrac{1}{120}\,F_{6,1,0}
-
\tfrac{1}{48}\,F_{6,0,0}
+
\tfrac{1}{72}\,F_{4,1,0}^2
\big)\,
x^6y^2
+
\big(
\tfrac{1}{48}\,
F_{5,1,0}
+
\tfrac{1}{120}\,
F_{6,0,1}
\big)\,
x^6yz
\\
&
\ \ \ \ \
+
{\rm O}_{y,z}(3)
+
{\rm O}_{x,y,z}(9).
\endaligned
\]
The same prenormalization holds for $v = G(r,s,t)$, of course.

According to~{\cite[Sec.~25]{Merker-2022}}, 
already at order 6, 
the stability group is 1-dimensional:
\[
\left[
\begin{array}{cccc}
a_{1,1} & \red{\bf 0} & \red{\bf 0} & \red{\bf 0}
\\
\red{\bf 0} & \red{\bf 1} & \red{\bf 0} & \red{\bf 0}
\\
\red{\bf 0} & \red{\bf 0} & \red{\tfrac{1}{a_{1,1}}} & \red{\bf 0}
\\
\red{\bf 0} & \red{\bf 0} & \red{\bf 0} & \red{a_{1,1}^2}
\end{array}
\right]^{\green{\bf 6}}.
\]

Therefore, $F_{4,1,0} \propto G_{4,1,0}$ is a relative invariant,
the lowest order one in fact, and all other
Taylor coefficients also are relative invariants, 
obviously. In fact:
\[
0
\overset{\green{\bf 4,1,0}}{\,\,=\,\,}
-\,\tfrac{1}{24}\,
F_{4,1,0}\,a_{1,1}^2
+
\tfrac{1}{24}\,
G_{4,1,0}\,a_{1,1}^4.
\]

Consequently, we must open two branches:
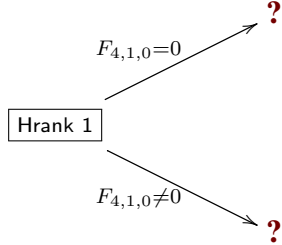
\[
\xymatrix{
&&
\red{\text{\bf ?}}
\\
\ar[urr]^{F_{4,1,0}=0}
\boxed{\scriptstyle\sf Hrank~1}
\ar[drr]_{F_{4,1,0}\neq0}
\\
&&
\red{\text{\bf ?}}
}
\]

\begin{Proposition}
In the branch $F_{4,1,0} \neq 0$, there are no affinely 
homogeneous models.
\end{Proposition}

Before starting the proof, 
without presenting the details,
let us state up to order 6, 
that $\eqL$ gives the following value:
\[
\aligned
L
&
\,=\,
\ \
\Big(
\big[
1-y-\tfrac{1}{2}\,F_{4,1,0}\,u
\big]\,T_1
+
u\,T_3
+
x\,A_{1,1}
\Big)\,\tfrac{\partial}{\partial x}
\\
&
\ \ \ \ \
+
\Big(
\big[
\tfrac{1}{2}\,F_{4,1,0}\,x
-
z
+
\tfrac{1}{5}\,F_{6,0,1}\,u
-
\tfrac{1}{5}\,F_{5,1,0}\,u
\big]\,T_1
+
\big[
1-y+\tfrac{1}{2}\,u\,F_{4,1,0}
\big]\,T_2
-
\tfrac{4}{3}\,x\,T_3
\Big)\,\tfrac{\partial}{\partial y}
\\
&
\ \ \ \ \
+
\Big(
\big[
-\tfrac{3}{10}\,F_{6,0,1}\,x
+
\tfrac{3}{10}\,F_{5,1,0}\,x
-
F_{4,1,0}\,y
-
\tfrac{1}{10}\,F_{6,0,0}\,u
-
\tfrac{1}{4}\,F_{4,1,0}^2\,u
\big]\,T_1
\\
&
\ \ \ \ \ \ \ \ \ \ \ \ \ \ \ \ \ \ \ \ \ \ \ 
+
\big[
-F_{4,1,0}\,x
-
2\,z
-
\tfrac{1}{10}\,F_{5,1,0}\,u
\big]\,T_2
+
\big[
1
-
y
+
\tfrac{2}{3}\,F_{4,1,0}\,u
\big]\,T_3
-
z\,A_{1,1}
\Big)\,\tfrac{\partial}{\partial z}
\\
&
\ \ \ \ \
+
\Big(
x\,T_1
+
u\,T_2
+
2\,u\,A_{1,1}
\Big)\,\tfrac{\partial}{\partial u},
\endaligned
\]
where $T_1, T_2, T_3$ and $A_{1,1}$ are free parameters.
 
Furthermore, at order 5, there remains 1 equation which
behaves differently in the two branches:
\[
0
\overset{\green{\bf 4,1,0}}{\,\,=\,\,}
\tfrac{1}{24}\,F_{5,1,0}\,T_1
+
\tfrac{1}{8}\,F_{4,1,0}\,T_2
+
\tfrac{1}{12}\,F_{4,1,0}\,A_{1,1},
\]
since $A_{1,1}$ may be solved only if $0 \neq F_{4,1,0}$, 
and there remain 2 equations at order 6:
\[
\aligned
0
&
\overset{\green{\bf 6,0,0}}{\,\,=\,\,}
\big(
\tfrac{1}{720}\,
F_{7,0,0}
+
\tfrac{1}{240}\,
F_{4,1,0}\,F_{6,0,1}
\big)\,T_1
+
\big(
\tfrac{1}{720}\,F_{6,1,0}
-
\tfrac{1}{720}\,F_{6,0,0}
+
\tfrac{1}{96}\,F_{4,1,0}^2
\big)\,T_2
\\
&
\ \ \ \ \
+
\big(
\tfrac{1}{720}\,F_{6,0,1}
-
\tfrac{1}{90}\,F_{5,1,0}
\big)\,T_3
+
\tfrac{1}{180}\,
F_{6,0,0}\,
A_{1,1},
\\
0
&
\overset{\green{\bf 5,1,0}}{\,\,=\,\,}
\big(
\tfrac{1}{120}\,F_{6,1,0}
-
\tfrac{1}{24}\,F_{6,0,0}
+
\tfrac{1}{48}\,F_{4,1,0}^2
\big)\,T_1
+
\tfrac{1}{30}\,F_{5,1,0}\,T_2
-
\tfrac{1}{72}\,F_{4,1,0}\,T_3
+
\tfrac{1}{40}\,F_{5,1,0}\,A_{1,1}.
\endaligned
\]

\proof
If $F_{4,1,0} \neq 0$, looking at 
$\overset{\green{\bf 4,1,0}}{\,=\,}$ of $\eqFG$
above, we can normalize:
\[
G_{4,1,0}
\,=\,
\pm\,1,
\]
and symetrically $F_{4,1,0} = \pm 1$. Then at the 
infinitesimal level, we may solve
from $\overset{\green{\bf 4,1,0}}{\,=\,}$
of $\eqL$:
\[
A_{1,1}
\,=\,
\mp\,\tfrac{1}{2}\,F_{5,1,0}\,T_1
-
\tfrac{3}{2}\,T_2,
\]
whence by replacement in
$\overset{\green{\bf 5,1,0}}{\,=\,}$ of $\eqL$:
\[
\aligned
0
&
\overset{\green{\bf 5,1,0}}{\,\,=\,\,}
\big(
\tfrac{1}{120}\,F_{6,1,0}
-
\tfrac{1}{24}\,F_{6,0,0}
\mp
\tfrac{1}{80}\,F_{5,1,0}^2
+
\tfrac{1}{48}
\big)\,T_1
-
\tfrac{1}{240}\,F_{5,1,0}\,T_2
\mp
\tfrac{1}{72}\,T_3.
\\
\endaligned
\]

This always is a contradictory
nontrivial linear relation between 
the transitivity parameters $\{T_1, T_2, T_3\}$,
because $\mp \frac{1}{72} \neq 0$.
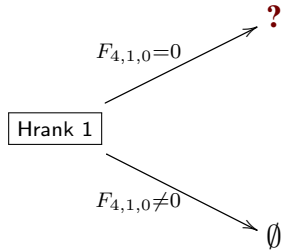
\endproof

\[
\xymatrix{
&&
\red{\text{\bf ?}}
\\
\ar[urr]^{F_{4,1,0}=0}
\boxed{\scriptstyle\sf Hrank~1}
\ar[drr]_{F_{4,1,0}\neq0}
\\
&&
\emptyset
}
\]

Therefore, $F_{4,1,0} = 0$ necessarily.
At order 6, $\eqL$ consists of 2 equations:
\[
\aligned
0
&
\overset{\green{\bf 6,0,0}}{\,\,=\,\,}
-\,\tfrac{1}{720}\,
F_{6,0,0}\,
a_{1,1}^2
+
\tfrac{1}{720}\,
G_{6,0,0}\,
a_{1,1}^6,
\\
0
&
\overset{\green{\bf 5,1,0}}{\,\,=\,\,}
-\,\tfrac{1}{120}\,
F_{5,1,0}\,
a_{1,1}^2
+
\tfrac{1}{120}\,
G_{5,1,0}\,
a_{1,1}^5.
\endaligned
\]
Again, we must open two branches.

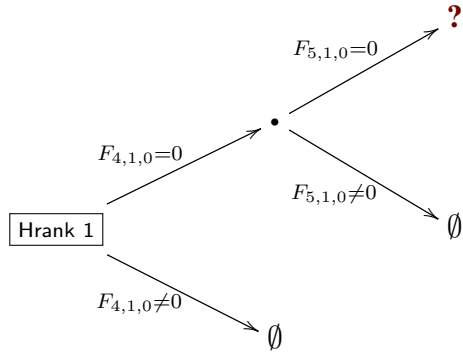
\[
\xymatrix{
&&
&&
\red{\text{\bf ?}}
\\
&&
\ar[urr]^{F_{5,1,0}=0}
\centersmallbullet
\ar[drr]_{F_{5,1,0}\neq0}
&&
\\
\ar[urr]^{F_{4,1,0}=0}
\boxed{\scriptstyle\sf Hrank~1}
\ar[drr]_{F_{4,1,0}\neq0}
&&
&&
\emptyset
\\
&&
\emptyset
&&
}
\]

But quickly, $\overset{\green{\bf 4,1,0}}{\,=\,}$ of $\eqL$ above 
\[
0
\overset{\green{\bf 4,1,0}}{\,\,=\,\,}
\tfrac{1}{24}\,F_{5,1,0}\,T_1
+
0
+
0,
\]
forces $F_{5,1,0} = 0$, so that one branch is void.

Similarly, the relative invariancy of $F_{6,0,0}$ creates 
two branches:
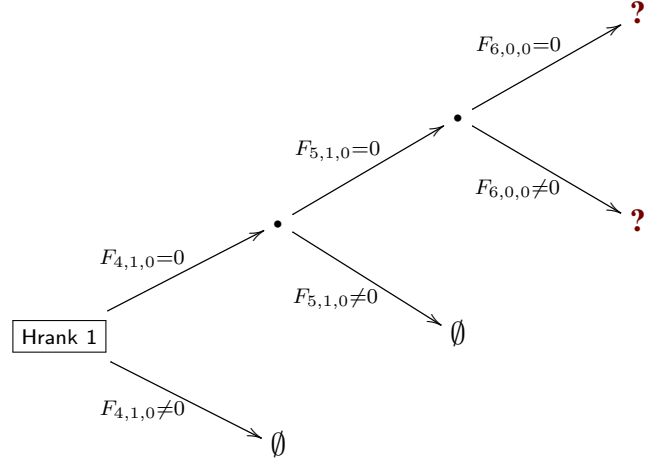
\[
\xymatrix{
&&
&&
&&
\red{\text{\bf ?}}
\\
&&
&&
\ar[urr]^{F_{6,0,0}=0}
\centersmallbullet
\ar[drr]_{F_{6,0,0}\neq0}
\\
&&
\ar[urr]^{F_{5,1,0}=0}
\centersmallbullet
\ar[drr]_{F_{5,1,0}\neq0}
&&
&&
\red{\text{\bf ?}}
\\
\ar[urr]^{F_{4,1,0}=0}
\boxed{\scriptstyle\sf Hrank~1}
\ar[drr]_{F_{4,1,0}\neq0}
&&
&&
\emptyset
\\
&&
\emptyset
&&
}
\]

Since $\overset{\green{\bf 5,1,0}}{\,=\,}$ of $\eqL$ becomes:
\[
0
\overset{\green{\bf 5,1,0}}{\,\,=\,\,}
\big(
\tfrac{1}{120}\,F_{6,1,0}
-
\tfrac{1}{24}\,F_{6,0,0}
+
0
\big)\,T_1
+
0
+
0
+
0,
\]
we have: 
\[
F_{6,1,0} 
\,=\,
5\,F_{6,0,0}.
\]

Then at orders 6 and 7, $\eqL$ consists of:
\[
\aligned
0
&
\overset{\green{\bf 6,0,0}}{\,\,=\,\,}
\tfrac{1}{720}\,F_{7,0,0}\,T_1
+
\tfrac{1}{180}\,F_{6,0,0}\,T_2
+
\tfrac{1}{720}\,F_{6,0,1}\,T_3
+
\tfrac{1}{180}\,F_{6,0,0}\,A_{1,1},
\\
0
&
\overset{\green{\bf 7,0,0}}{\,\,=\,\,}
\big(
\tfrac{1}{5040}\,F_{8,0,0}
-
\tfrac{1}{2400}\,F_{6,0,1}^2
\big)\,T_1
+
\big(
-\tfrac{1}{5040}\,F_{7,0,0}
+
\tfrac{1}{5040}\,F_{7,1,0}
\big)\,T_2
\\
&
\ \ \ \ \ \ \ \ \ \ \ \ \ \ \ \ \ \ \ \ \ \ \ \ \ \ \ \ \ \ \ \ \ \ \
\ \ \ \ \ \ \ \ \ \ \ \ \ \ \ \ \ \ \ \ 
+
\big(
-\tfrac{1}{270}\,F_{6,0,0}
+
\tfrac{1}{5040}\,F_{7,0,1}
\big)\,T_3
+
\tfrac{1}{1008}\,F_{7,0,0}\,A_{1,1},
\\
0
&
\overset{\green{\bf 6,1,0}}{\,\,=\,\,}
\big(
-\tfrac{1}{720}\,F_{7,0,0}
+
\tfrac{1}{720}\,F_{7,1,0}
\big)\,T_1
+
\tfrac{1}{36}\,F_{6,0,0}\,T_2
+
\tfrac{1}{144}\,F_{6,0,1}\,T_3
+
\tfrac{1}{36}\,F_{6,0,0}\,A_{1,1},
\\
0
&
\overset{\green{\bf 6,0,1}}{\,\,=\,\,}
\big(
\tfrac{1}{720}\,F_{7,0,1}
-
\tfrac{1}{45}\,F_{6,0,0}
\big)\,T_1
+
\tfrac{1}{240}\,F_{6,0,1}\,T_2
+
\tfrac{1}{240}\,F_{6,0,1}\,A_{1,1}.
\endaligned
\]

\begin{Proposition}
In the branch $F_{6,0,0} \neq 0$, there are no affinely homogenous
models.
\end{Proposition}

\proof
From $\overset{\green{\bf 6,0,0}}{\,=\,}$ of $\eqFG$, 
we may normalize:
\[
G_{6,0,0}
\,=\,
\pm1
\,=\,
F_{6,0,0},
\]
then from $\overset{\green{\bf 6,0,0}}{\,=\,}$ of $\eqL$:
\[
A_{1,1}
\,:=\,
\mp
\tfrac{1}{4}\,F_{7,0,0}\,T_1
-
T_2
\mp
\tfrac{1}{4}\,F_{6,0,1}\,T_3.
\]

A replacement gives:
\[
0
\overset{\green{\bf 6,0,1}}{\,\,=\,\,}
\mp\,
\tfrac{1}{960}\,F_{6,0,1}^2\,T_3
+
\big(
\tfrac{1}{720}\,F_{7,0,1}
\mp
\tfrac{1}{960}\,F_{7,0,0}\,F_{6,0,1}
\mp
\tfrac{1}{45}
\big)\,T_1,
\]
whence $F_{6,0,1} = 0$ necessarily, and then:

\[
\aligned
0
&
\overset{\green{\bf 7,0,0}}{\,\,=\,\,}
\ast\,T_1
+
\ast\,T_2
+
\big(
\tfrac{1}{5040}\,F_{7,0,1}
\mp
\tfrac{1}{270}
\big)\,T_3,
\\
0
&
\overset{\green{\bf 6,0,1}}{\,\,=\,\,}
0
+
\big(
\tfrac{1}{720}\,F_{7,0,1}
\mp
\tfrac{1}{45}
\big)\,T_1,
\endaligned
\]
where $\ast$ are unimportant, but this gives 
the two
noncoinciding values 
$\pm \frac{56}{3}$ and $\pm 16$ 
for $F_{7,0,1}$. 
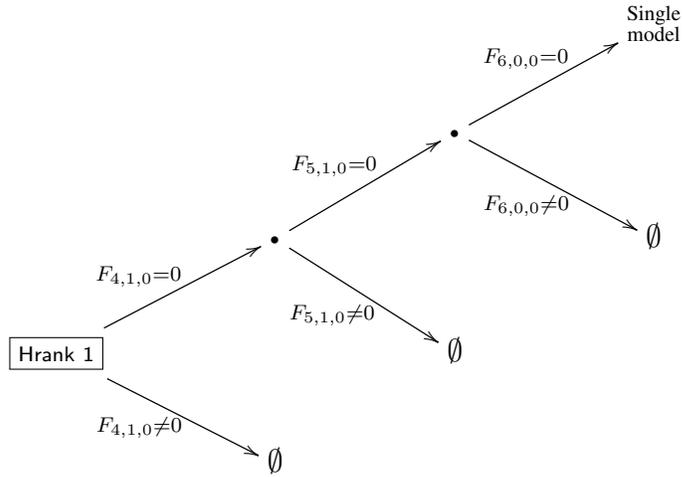
\endproof

\[
\xymatrix{
&&
&&
&&
{\substack{
\text{Single}
\\
\text{model}}}
\\
&&
&&
\ar[urr]^{F_{6,0,0}=0}
\centersmallbullet
\ar[drr]_{F_{6,0,0}\neq0}
\\
&&
\ar[urr]^{F_{5,1,0}=0}
\centersmallbullet
\ar[drr]_{F_{5,1,0}\neq0}
&&
&&
\emptyset
\\
\ar[urr]^{F_{4,1,0}=0}
\boxed{\scriptstyle\sf Hrank~1}
\ar[drr]_{F_{4,1,0}\neq0}
&&
&&
\emptyset
\\
&&
\emptyset
&&
}
\]

\begin{Theorem}
\label{Thm-H3-R4}
Among constant Hessian rank 1 hypersurfaces
$H^3 \subset \R^4$, there is 
a single affinely homogeneous model, 
lying in the branch $F_{2,1} \neq 0$, $F_{3,1} = 0$, $F_{5,0} = 0$,
of equation:
\[
\footnotesize
\aligned
u
&
\,=\,
\ \
\tfrac{x^2}{2}
\\
&
\ \ \ \ \
+
\tfrac{x^2y}{2}
\\
&
\ \ \ \ \
+
\tfrac{x^3z}{6}
+
\tfrac{x^2y^2}{2}
\\
&
\ \ \ \ \
+
\tfrac{x^3yz}{2}
+
\tfrac{x^2y^3}{2}
\\
&
\ \ \ \ \
+
\tfrac{1}{8}\,x^4z^2
+
x^3y^2z
+
\tfrac{1}{2}\,x^2y^4
\\
&
\ \ \ \ \
+
\tfrac{5}{8}\,x^4yz^2
+
\tfrac{5}{3}\,x^3y^3z
+
\tfrac{1}{2}\,x^2y^5,
\\
&
\ \ \ \ \
+
\tfrac{1}{8}\,x^5z^3
+
\tfrac{15}{8}\,x^4y^2z^2
+
\tfrac{5}{2}\,x^3y^4z
+
\tfrac{1}{2}\,x^2y^6
\\
&
\ \ \ \ \
+
\tfrac{7}{8}\,x^5yz^3
+
\tfrac{35}{8}\,x^4y^3z^2
+
\tfrac{7}{2}\,x^3y^5z
+
\tfrac{1}{2}\,x^2y^7
\\
&
\ \ \ \ \
+
\tfrac{7}{48}\,x^6z^4
+
\tfrac{7}{2}\,x^5y^2z^3
+
\tfrac{35}{4}\,x^4y^4z^2
+
\tfrac{14}{3}\,x^3y^6z
+
\tfrac{1}{2}\,x^2y^8
+
\\
&
\ \ \ \ \
+
{\rm O}_{x,y,z}(11),
\endaligned
\]
with 4-dimensional affine symmetry algebra generated by:
\[
\aligned
e_1
&
\,:=\,
(1-y)\,\partial_x
-
z\,\partial_y
+
x\,\partial_u,
\\
e_2
&
\,:=\,
(1-y)\,\partial_y
-
2z\,\partial_z
+
u\,\partial_u,
\\
e_3
&
\,:=\,
u\,\partial_x
-
\tfrac{4}{3}\,x\,\partial_y
+
(1-y)\,\partial_z,
\\
e_4
&
\,:=\,
x\,\partial_x
-
z\,\partial_z
+
2\,u\,\partial_u.
\endaligned
\]
\end{Theorem}

\proof
Putting $F_{6,0,0} := 0$, 
and knowing $F_{6,1,0} = 5\, F_{6,0,0} = 0$,
at order 6 for $\eqL$, it remains only:
\[
0
\overset{\green{\bf 6,0,0}}{\,\,=\,\,}
\tfrac{1}{720}\,
F_{7,0,0}\,T_1
+
\tfrac{1}{720}\,
F_{6,0,1}\,T_3,
\]
whence $0 = F_{7,0,0} = F_{6,0,1}$.

At order 7, $\eqL$ reads:
\[
\aligned
0
&
\overset{\green{\bf 7,0,0}}{\,\,=\,\,}
\tfrac{1}{5040}\,
F_{8,0,0}\,T_1
+
\tfrac{1}{5040}\,F_{7,1,0}\,T_2
+
\tfrac{1}{5040}\,F_{7,0,1}\,T_3,
\\
0
&
\overset{\green{\bf 6,1,0}}{\,\,=\,\,}
\tfrac{1}{720}\,F_{7,1,0}\,T_1,
\\
0
&
\overset{\green{\bf 6,0,1}}{\,\,=\,\,}
\tfrac{1}{720}\,F_{7,0,1}\,T_1,
\endaligned
\]
whence $F_{8,0,0} = F_{7,1,0} = F_{7,0,1}$.

Generally, one can see that for all $\mu \geqslant 7$:
\[
0
\,=\,
F_{\mu,0,0}
\,=\,
F_{\mu-1,1,0}
\,=\,
F_{\mu,0,1}.
\qedhere
\]
\endproof

\begin{Corollary}
A closed expression for the graphing function 
$F(x,y,z)$ is:
\[
u
\,=\,
\frac{1}{3\,z^2}
\Big\{
\big(
1-2\,y+y^2-2\,xz
\big)^{3/2}
-
(1-y)\,
\big(
1-2\,y+y^2-3\,xz
\big)
\Big\}.
\]
\end{Corollary}

\proof
By expanding the numerator in power series, one realizes
that the singularity $\frac{1}{z^2}$ is removable,
and that the power series
expansion matches with that
of Theorem~{\ref{Thm-H3-R4}} up to order 10 monomials.

On the other hand, one verifies that 
$e_1$, $e_2$, $e_3$, $e_4$ are infinitesimal
symmetries of this closed form.
\endproof

\SectionHead{Fourfolds $H^4 \subset \R^5$}
{fourfolds-H4-R5}

In $\R^5$, consider an affine-linear map
$(x,y,z,w,u) \longmapsto (r,s,t,p,v)$ fixing the origin:
\[
\aligned
r
&
\,:=\,
a_{1,1}\,x+a_{1,2}\,y+a_{1,3}\,z+a_{1,4}\,w+b_1\,u,
\\
s
&
\,:=\,
a_{2,1}\,x+a_{2,2}\,y+a_{2,3}\,z+a_{2,4}\,w+b_2\,u,
\\
t
&
\,:=\,
a_{3,1}\,x+a_{3,2}\,y+a_{3,3}\,z+a_{3,4}\,w+b_3\,u,
\\
p
&
\,:=\,
a_{4,1}\,x+a_{4,2}\,y+a_{4,3}\,z+a_{4,4}\,w+b_4\,u,
\\
v
&
\,:=\,
c_1\,x+c_2\,y+c_3\,z+c_4\,w+d\,u,
\endaligned
\ \ \ \ \
\text{with}
\ \ \ \ \
0
\,\neq\,
\left\vert\!
\begin{array}{ccccc}
a_{1,1} & a_{1,2} & a_{1,3} & a_{1,4} & b_1
\\
a_{2,1} & b_{2,2} & a_{2,3} & a_{2,4} & b_2
\\
a_{3,1} & b_{3,2} & a_{3,3} & a_{3,4} & b_3
\\
a_{4,1} & b_{4,2} & a_{4,3} & a_{4,4} & b_4
\\
c_1 & c_2 & c_3 & c_4 & d
\end{array}
\!\right\vert.
\]
Also, consider two graphed analytic hypersurfaces:
\[
u
\,=\,
F(x,y,z,w)
\ \ \ \ \ \ \ \ 
{\scriptstyle{(F(0,0,0,0)\,=\,0)}}
\ \ \ \ \ \ \ \ \ \ \ \
\text{and}
\ \ \ \ \ \ \ \ \ \ \ \ 
v
\,=\,
G(r,s,t,p)
\ \ \ \ \ \ \ \ 
{\scriptstyle{(0\,=\,G(0,0,0,0))}},
\]
and assume that the above map is an affine equivalence 
$\{ u = F\} \longrightarrow \{ v = G\}$.

The main hypothesis of
constant Hessian rank 1, 
after elementary preliminary
transformations:
\[
u
\,=\,
\tfrac{x^2}{2}
+
{\rm O}_{x,y,z,w}(3),
\]
reads as:
\[
\aligned
1
\,\equiv\,
\rank\,
\left[\!
\begin{array}{cccc}
F_{xx} & F_{xy} & F_{xz} & F_{xw}
\\
F_{yx} & F_{yy} & F_{yz} & F_{yw}
\\
F_{zx} & F_{zy} & F_{zz} & F_{zw}
\\
F_{wx} & F_{wy} & F_{wz} & F_{ww}
\end{array}
\!\right],
\endaligned
\]
which is then equivalent to:
\[
\aligned
0
&
\,\equiv\,
\left\vert\!
\begin{array}{cc}
F_{xx} & F_{xy} 
\\
F_{yx} & F_{yy}
\end{array}
\!\right\vert
\,\equiv\,
\left\vert\!
\begin{array}{cc}
F_{xx} & F_{xz} 
\\
F_{yx} & F_{yz}
\end{array}
\!\right\vert
\,\equiv\,
\left\vert\!
\begin{array}{cc}
F_{xx} & F_{xw} 
\\
F_{yx} & F_{yw}
\end{array}
\!\right\vert
\\
&
\,\equiv\,
\left\vert\!
\begin{array}{cc}
F_{xx} & F_{xy} 
\\
F_{zx} & F_{zy}
\end{array}
\!\right\vert
\,\equiv\,
\left\vert\!
\begin{array}{cc}
F_{xx} & F_{xz} 
\\
F_{zx} & F_{zz}
\end{array}
\!\right\vert
\,\equiv\,
\left\vert\!
\begin{array}{cc}
F_{xx} & F_{xw} 
\\
F_{zx} & F_{zw}
\end{array}
\!\right\vert,
\\
&
\,\equiv\,
\left\vert\!
\begin{array}{cc}
F_{xx} & F_{xy} 
\\
F_{wx} & F_{wy}
\end{array}
\!\right\vert
\,\equiv\,
\left\vert\!
\begin{array}{cc}
F_{xx} & F_{xz} 
\\
F_{wx} & F_{wz}
\end{array}
\!\right\vert
\,\equiv\,
\left\vert\!
\begin{array}{cc}
F_{xx} & F_{xw} 
\\
F_{wx} & F_{ww}
\end{array}
\!\right\vert.
\endaligned
\]
By affine invariancy of the Hessian matrix rank, 
the same holds about $v = \tfrac{r^2}{2} + {\rm O}_{r,s,t}(3)$.

The fundamental equation which 
holds identically in $\R\{x,y,z,w\}$:
\[
\aligned
0
\,\equiv\,
\eqFG(x,y,z,w),
\endaligned
\]
writes:
\[
\aligned
\eqFG
&
\,:=\,
-\,c_1\,x-c_2\,y-c_3\,z-c_4\,w-d\,F(x,y,z)
\\
&
\ \ \ \ \
+
G
\Big(
a_{1,1}x+a_{1,2}y+a_{1,3}z+a_{1,4}w+b_1F(x,y,z,w),\,\,
\\
&
\ \ \ \ \ \ \ \ \ \ \ \ \ \ \ \ 
a_{2,1}x+a_{2,2}y+a_{2,3}z+a_{2,4}w+b_2F(x,y,z,w),\,\,
\\
&
\ \ \ \ \ \ \ \ \ \ \ \ \ \ \ \ 
a_{3,1}x+a_{3,2}y+a_{3,3}z+a_{3,4}w+b_3F(x,y,z,w),\,\,
\\
&
\ \ \ \ \ \ \ \ \ \ \ \ \ \ \ \ 
a_{4,1}x+a_{4,2}y+a_{4,3}z+a_{4,4}w+b_4F(x,y,z,w)
\Big).
\endaligned
\]

Also, an affine vector field:
\[
\aligned
L
&
\,=\,
\ \ 
\big(
T_1+A_{1,1}\,x+A_{1,2}\,y+A_{1,3}\,z+A_{1,4}\,w+B_1\,u
\big)\,\frac{\partial}{\partial x}
\\
&
\ \ \ \ \
+
\big(
T_2+A_{2,1}\,x+A_{2,2}\,y+A_{2,3}\,z+A_{2,4}\,w+B_2\,u
\big)\,\frac{\partial}{\partial y}
\\
&
\ \ \ \ \
+
\big(
T_3+A_{3,1}\,x+A_{3,2}\,y+A_{3,3}\,z+A_{3,4}\,w+B_3\,u
\big)\,\frac{\partial}{\partial z}
\\
&
\ \ \ \ \
+
\big(
T_4+A_{4,1}\,x+A_{4,2}\,y+A_{4,3}\,z+A_{4,4}\,w+B_4\,u
\big)\,\frac{\partial}{\partial w}
\\
&
\ \ \ \ \
+
\big(
T_0+C_1\,x+C_2\,y+C_3\,z+C_4\,w+D\,u
\big)\,\frac{\partial}{\partial u},
\endaligned
\]
is tangent to $\{u = F(x,y,z,w)\}$ if and only if:
\[
\aligned
0
&
\,\equiv\,
\eqL(x,y,z,w)
\\
&
\,=:\,
L
\big(
-\,u+F(x,y,z,w)
\big)
\Big\vert_{u=F(x,y,z,w)},
\endaligned
\]
identically as power series in $\R\{x,y,z,w\}$.

According to Theorems~1.4, 13.1, 1.5, 25.2
in~{\cite{Merker-2022}}, 
if $H^4 \subset \R^5$ is not affinely equivalent to a product
with $\R^1$ or $\R^2$ or $\R^3$, 
its graphing function $F(x,y,z,w)$ can be 
{\sl pre-normalized}\,\,---\,\,that is, 
{\em normalized before creating any branching}\,\,---\,\,up 
to order $4 + 5 = 9$ included 
and modulo ${\rm O}_{y,z,w}(3)$ as:
\[
\footnotesize
\aligned
u
&
\,=\,
\ \
\tfrac{x^2}{2}
\\
&
\ \ \ \ \
+
\tfrac{x^2y}{2}
\\
&
\ \ \ \ \
+
\tfrac{x^3z}{6}
+
\tfrac{x^2y^2}{2}
\\
&
\ \ \ \ \
+
\tfrac{x^4w}{24}
+
\tfrac{x^3yz}{2}
+
\tfrac{x^2y^3}{2}
\\
&
\ \ \ \ \
+
F_{5,1,0,0}\,\tfrac{x^5y}{120}
+
\tfrac{x^4yw}{6}
+
\tfrac{x^4z^2}{8}
+
x^3y^2z
\\
&
\ \ \ \ \
+
F_{7,0,0,0}\,\tfrac{x^7}{5040}
+
F_{6,1,0,0}\,\tfrac{x^6y}{720}
+
F_{6,0,0,1}\,\tfrac{x^6w}{720}
+
F_{5,1,0,0}\,\tfrac{x^5y^2}{24}
+
\tfrac{x^5zw}{12}
\\
&
\ \ \ \ \
+
\tfrac{1}{40320}\,F_{8,0,0,0}\,x^8
+
\tfrac{1}{5040}\,F_{7,1,0,0}\,x^7y
+
\tfrac{1}{5040}\,F_{7,0,1,0}\,x^7z
+
\tfrac{1}{5040}\,F_{7,0,0,1}\,x^7w
\\
&
\ \ \ \ \ \ \ \ \ \ \ \ \ \ \ \ \ \ \ \ \ \ \ \ \ \ \ \ \ \ \ \ \ \ \
\ \ \ \ \
+
\tfrac{1}{120}\,F_{6,1,0,0}\,x^6y^2
+
\tfrac{1}{48}\,F_{5,1,0,0}\,x^6yz
+
\tfrac{1}{120}\,F_{6,0,0,1}\,x^6yw
+
\tfrac{1}{72}\,x^6w^2
\\
&
\ \ \ \ \
+
\tfrac{1}{362880}\,F_{9,0,0,0}\,x^9
+
\tfrac{1}{40320}\,F_{8,1,0,0}\,x^8y
+
\tfrac{1}{40320}\,F_{8,0,1,0}\,x^8z
+
\tfrac{1}{40320}\,F_{8,0,0,1}\,x^8w
\\
&
\ \ \ \ \ \ \ \ \ \ \ \ \ \ \ \ \ \ \ \ \ \ \ \ \ \ \ \ \ \ \ \ \ \ \
\ \ \ \ \
+
\tfrac{1}{10080}
\big(
14\,F_{7,1,0,0}
-
42\,F_{7,0,0,0}
\big)\,
x^7y^2
+
\tfrac{1}{5040}\,
\big(
7\,F_{7,0,1,0}
+
21\,F_{6,1,0,0}
\big)\,x^7yz
\\
&
\ \ \ \ \ \ \ \ \ \ \ \ \ \ \ \ \ \ \ \ \ \ \ \ \ \ \ \ \ \ \ \ \ \ \
\ \ \ \ \
+
\tfrac{1}{5040}\,
\big(
7\,F_{7,0,0,1}
+
35\,F_{5,1,0,0}
\big)\,
x^7yw
+
\tfrac{1}{240}\,F_{6,0,0,1}\,x^7zw
\\
&
\ \ \ \ \
+
{\rm O}_{y,z,w}(3)
+
{\rm O}_{x,y,z,w}(10).
\endaligned
\]

According to~{\cite[Sec.~25]{Merker-2022}}, 
already at order 7, 
the stability group is 1-dimensional:
\[
\left[
\begin{array}{ccccc}
a_{1,1} & \red{\bf 0} & \red{\bf 0} & \red{\bf 0} & \red{\bf 0}
\\
\red{\bf 0} & \red{\bf 1} & \red{\bf 0} & \red{\bf 0} & \red{\bf 0}
\\
\red{\bf 0} & \red{\bf 0} & \red{\tfrac{1}{a_{1,1}}} & \red{\bf 0} & 
\red{\bf 0}
\\
\red{\bf 0} & \red{\bf 0} & \red{\bf 0} & 
\red{\tfrac{1}{a_{1,1}^2}} & \red{\bf 0}
\\
\red{\bf 0} & \red{\bf 0} & \red{\bf 0} & \red{\bf 0} & \red{a_{1,1}^2}
\end{array}
\right]^{\green{\bf 7}}.
\]

Moreover, $F_{5,1,0,0} \propto G_{5,1,0,0}$ is a relative invariant,
the lowest order one in fact, and all other
Taylor coefficients also are relative invariants, 
obviously. In fact:
\[
0
\overset{\green{\bf 5100}}{\,\,=\,\,}
-\,\tfrac{1}{120}\,
F_{5,1,0,0}\,a_{1,1}^2
+
\tfrac{1}{120}\,
G_{5,1,0,0}\,a_{1,1}^5.
\]
For the moment, we do not open a branching here.

Up to order 6, $\eqL$ gives:
\[
\footnotesize
\aligned
L
&
\,=\,
\ \
\Big(
T_1
+
A_{1,1}\,x
-
T_1\,y
+
\big[
\tfrac{1}{5}\,F_{6,0,0,1}\,T_1
+
\tfrac{2}{3}\,T_3
\big]\,
u
\Big)\,
\tfrac{\partial}{\partial x}
\\
&
\ \ \ \ \ 
+
\Big(
T_2
+
\big[
-T_3
-
\tfrac{1}{5}\,F_{6,0,0,1}\,T_1
\big]\,x
-
T_2\,y
-
T_1\,z
+
\big[
\tfrac{1}{2}\,T_4
-
\tfrac{1}{5}\,F_{5,1,0,0}\,T_1
\big]\,u
\Big)\,
\tfrac{\partial}{\partial y}
\\
&
\ \ \ \ \ 
+
\Big(
T_3
+
\big[
\tfrac{3}{10}\,F_{5,1,0,0}\,T_1
-
T_4
\big]\,x
-
T_3\,y
+
\big[
-2\,T_2
-
A_{1,1}
\big]\,z
-
T_1\,w
\\
&
\ \ \ \ \ \ \ \ \ \ \ \ \ \ \ \ \ \ \ \ \ \ \ \ \ \ \ \ \ \ \ \ \ \ \
+
\big[
\big(
\tfrac{1}{10}\,F_{7,0,1,0}
-
\tfrac{1}{10}\,F_{6,1,0,0}
+
\tfrac{1}{25}\,F_{6,0,0,1}^2
\big)\,T_1
+
\tfrac{3}{10}\,F_{5,1,0,0}\,T_2
-
\tfrac{1}{15}\,F_{6,0,0,1}\,T_3
\big]\,u
\Big)\,
\tfrac{\partial}{\partial z}
\\
&
\ \ \ \ \ 
+
\Big(
T_4
+
\big[
\big(
\tfrac{1}{5}\,F_{6,1,0,0}
-
\tfrac{1}{5}\,F_{7,0,1,0}
-
\tfrac{2}{25}\,F_{6,0,0,1}^2
\big)\,T_1
-
\tfrac{4}{5}\,F_{5,1,0,0}\,T_2
+
\tfrac{2}{15}\,F_{6,0,0,1}\,T_3
\big]\,x
+
\big[
-F_{5,1,0,0}\,T_1
-
T_4
\big]\,y
\\
&
\ \ \ \ \ \ \ \ \ \ \ \ \ \ \ \ \ \ \ \ \ \ \ \ \ \ \ \ \ \ \ \ \ \ \
+
\big[
\tfrac{2}{5}\,F_{6,0,0,1}\,T_1
-
\tfrac{2}{3}\,T_3
\big]\,z
+
\big[
-3\,T_2
-
2\,A_{1,1}
\big]\,w
\\
&
\ \ \ \ \ \ \ \ \ \ \ \ \ \ \ \ \ \ \ \ \ \ \ \ \ \ \ \ \ \ \ \ \ \ \
+
\big[
\big(
\tfrac{2}{25}\,F_{5,1,0,0}\,F_{6,0,0,1}
-
\tfrac{1}{15}\,F_{7,0,0,0}
\big)\,T_1
-
\tfrac{1}{15}\,F_{6,1,0,0}\,T_2
+
\tfrac{2}{5}\,F_{5,1,0,0}\,T_3
-
\tfrac{1}{15}\,F_{6,0,0,1}\,T_4
\big]\,u
\Big)\,
\tfrac{\partial}{\partial w}
\\
&
\ \ \ \ \ 
+
\Big(
T_1\,x
+
\big[
T_2+2\,A_{1,1}
\big]\,
u
\Big)\,
\tfrac{\partial}{\partial u},
\endaligned
\]
with the four transitivity parameters $T_1, T_2, T_3, T_4$,
plus a single possible isotropy parameter $A_{1,1}$.

Up to order 8, the remaining equations of $\eqL$ are:
\[
0
\overset{\green{\bf 5100}}{\,\,=\,\,}
\tfrac{1}{120}\,F_{6,1,0,0}\,T_1
+
\tfrac{1}{30}\,F_{5,1,0,0}\,T_2
+
\tfrac{1}{40}\,F_{5,1,0,0}\,A_{1,1},
\]
\[
\footnotesize
\aligned
0
&
\overset{\green{\bf 7000}}{\,\,=\,\,}
\Big(
\tfrac{1}{5040}\,
F_{8,0,0,0}
-
\tfrac{1}{3600}\,
F_{7,0,1,0}\,F_{6,0,0,1}
-
\tfrac{1}{9000}\,
F_{6,0,0,1}^3
-
\tfrac{1}{1200}\,
F_{5,1,0,0}^2
\Big)\,
T_1
\\
&
\ \ \ \ \ 
+
\Big(
\tfrac{1}{5040}\,F_{7,1,0,0}
-
\tfrac{1}{5040}\,F_{7,0,0,0}
-
\tfrac{1}{900}\,F_{6,0,0,1}\,F_{5,1,0,0}
\Big)\,
T_2
\\
&
\ \ \ \ \ 
+
\Big(
\tfrac{1}{5040}\,F_{7,0,1,0}
-
\tfrac{1}{720}\,F_{6,1,0,0}
+
\tfrac{1}{5040}\,F_{6,0,0,1}^2
\Big)\,T_3
+
\Big(
\tfrac{1}{5040}\,F_{7,0,0,1}
+
\tfrac{1}{480}\,F_{5,1,0,0}
\Big)\,T_4
+
\tfrac{1}{1008}\,F_{7,0,0,0}\,A_{1,1},
\endaligned
\]
\[
\footnotesize
\aligned
0
&
\overset{\green{\bf 6100}}{\,\,=\,\,}
\Big(
\tfrac{1}{720}\,F_{7,1,0,0}
-
\tfrac{1}{120}\,F_{7,0,0,0}
-
\tfrac{7}{1800}\,F_{5,1,0,0}\,F_{6,0,0,1}
\Big)\,T_1
+
\tfrac{1}{144}\,F_{6,1,0,0}\,T_2
-
\tfrac{1}{720}\,F_{5,1,0,0}\,T_3
+
\tfrac{1}{180}\,F_{6,1,0,0}\,A_{1,1}
\endaligned
\]
\[
\footnotesize
\aligned
0
&
\overset{\green{\bf 6001}}{\,\,=\,\,}
\Big(
\tfrac{1}{720}\,F_{7,0,0,1}
+
\tfrac{1}{240}\,F_{5,1,0,0}
\Big)\,T_1
+
\tfrac{1}{360}\,F_{6,0,0,1}\,T_2
-
\boxed{
\tfrac{1}{288}\,T_4}
+
\tfrac{1}{360}\,F_{6,0,0,1}\,A_{1,1},
\endaligned
\]
\[
\footnotesize
\aligned
0
&
\overset{\green{\bf 8000}}{\,\,=\,\,}
\Big(
\tfrac{1}{40320}\,F_{9,0,0,0}
-
\tfrac{1}{25200}\,
F_{7,0,1,0}\,F_{7,0,0,1}
-
\tfrac{1}{25200}\,F_{7,1,0,0}\,F_{6,0,0,1}
+
\tfrac{1}{16800}\,F_{7,0,1,0}\,F_{5,1,0,0}
+
\tfrac{1}{25200}\,F_{7,0,0,1}\,F_{6,1,0,0}
\\
&
\ \ \ \ \ \ \ \ \ \ \ \ \ \ \ \ \ \ \ \ \ \ \ \ \ \ \ \ \ \ \ \ \ \ \
-
\tfrac{1}{63000}\,F_{7,0,0,1}\,F_{6,0,0,1}^2
+
\tfrac{1}{7560}\,F_{7,0,0,0}\,F_{6,0,0,1}
+
\tfrac{1}{18000}\,F_{5,1,0,0}\,F_{6,0,0,1}^2
-
\tfrac{1}{7200}\,F_{5,1,0,0}\,F_{6,1,0,0}
\Big)\,T_1
\\
&
\ \ \ \ \
+
\Big(
\tfrac{1}{40320}\,F_{8,1,0,0}
-
\tfrac{1}{40320}\,F_{8,0,0,0}
-
\tfrac{1}{6300}\,F_{5,1,0,0}\,F_{7,0,0,1}
-
\tfrac{1}{21600}\,F_{6,1,0,0}\,F_{6,0,0,1}
\Big)\,T_2
\\
&
\ \ \ \ \
+
\Big(
\tfrac{1}{40320}\,F_{8,0,1,0}
-
\tfrac{1}{5040}\,F_{7,1,0,0}
+
\tfrac{1}{1680}\,F_{7,0,0,0}
+
\tfrac{1}{37800}\,F_{6,0,0,1}\,F_{7,0,0,1}
+
\tfrac{1}{3600}\,F_{5,1,0,0}\,F_{6,0,0,1}
\Big)\,T_3
\\
&
\ \ \ \ \
+
\Big(
\tfrac{1}{40320}\,F_{8,0,0,1}
-
\tfrac{1}{5040}\,F_{7,0,1,0}
+
\tfrac{1}{2880}\,F_{6,1,0,0}
-
\tfrac{1}{21600}\,F_{6,0,0,1}^2
\Big)\,T_4
+
\tfrac{1}{6720}\,F_{8,0,0,0}\,A_{1,1},
\endaligned
\]
\[
\footnotesize
\aligned
0
&
\overset{\green{\bf 7100}}{\,\,=\,\,}
\Big(
\tfrac{1}{5040}\,F_{8,1,0,0}
-
\tfrac{1}{5040}\,F_{8,0,0,0}
-
\tfrac{1}{5040}\,F_{5,1,0,0}\,F_{7,0,0,1}
-
\tfrac{1}{600}\,F_{7,0,1,0}\,F_{6,0,0,1}
\\
&
\ \ \ \ \ \ \ \ \ \ \ \ \ \ \ \ \ \ \ \ \ \ \ \ \ \ \ \ \ \ \ \ \ \ \
-
\tfrac{1}{1800}\,F_{6,1,0,0}\,F_{6,0,0,1}
-
\tfrac{3}{800}\,F_{5,1,0,0}^2
-
\tfrac{1}{1500}\,F_{6,0,0,1}^3
\Big)\,T_1
+
\\
&
\ \ \ \ \
+
\Big(
\tfrac{1}{420}\,F_{7,1,0,0}
-
\tfrac{1}{120}\,F_{7,0,0,0}
-
\tfrac{1}{150}\,F_{5,1,0,0}\,F_{6,0,0,1}
\Big)\,T_2
+
\Big(
\tfrac{1}{840}\,F_{7,0,1,0}
-
\tfrac{19}{2160}\,F_{6,1,0,0}
+
\tfrac{1}{900}\,F_{6,0,0,1}^2
\Big)\,T_3
\\
&
\ \ \ \ \
+
\Big(
\tfrac{1}{840}\,F_{7,0,0,1}
+
\tfrac{1}{90}\,F_{5,1,0,0}
\Big)\,T_4
+
\tfrac{1}{1008}\,F_{7,1,0,0}\,A_{1,1},
\endaligned
\]
\[
\footnotesize
\aligned
0
&
\overset{\green{\bf 7010}}{\,\,=\,\,}
\Big(
\tfrac{1}{5040}\,F_{8,0,1,0}
-
\tfrac{1}{5040}\,F_{7,1,0,0}
-
\tfrac{7}{2160}\,F_{7,0,0,0}
+
\tfrac{1}{12600}\,F_{6,0,0,1}\,F_{7,0,0,1}
-
\tfrac{1}{3600}\,F_{5,1,0,0}\,F_{6,0,0,1}
\Big)\,T_1
\\
&
\ \ \ \ \
+
\Big(
\tfrac{1}{1260}\,F_{7,0,1,0}
+
\tfrac{1}{1080}\,F_{6,1,0,0}
\Big)\,T_2
+
\Big(
-\tfrac{1}{7560}\,F_{7,0,0,1}
-
\tfrac{1}{720}\,F_{5,1,0,0}
\Big)\,T_3
+
\tfrac{1}{1080}\,F_{6,0,0,1}\,T_4
+
\tfrac{1}{1260}\,F_{7,0,1,0}\,A_{1,1},
\endaligned
\]
\[
\footnotesize
\aligned
0
&
\overset{\green{\bf 7001}}{\,\,=\,\,}
\Big(
\tfrac{1}{5040}\,F_{8,0,0,1}
-
\tfrac{1}{1120}\,F_{7,0,1,0}
+
\tfrac{1}{1440}\,F_{6,1,0,0}
-
\tfrac{1}{1200}\,F_{6,0,0,1}^2
\Big)\,T_1
\\
&
\ \ \ \ \
+
\Big(
\tfrac{1}{1680}\,F_{7,0,0,1}
-
\tfrac{1}{1440}\,F_{5,1,0,0}
\Big)\,T_2
+
\tfrac{1}{1680}\,F_{7,0,0,1}\,A_{1,1}.
\endaligned
\]

\begin{Observation}
$F_{6,0,0,1} \neq 0$, necessarily.
\end{Observation}

\proof
If we would have $F_{6,0,0,1} = 0$, 
because of the presence of $-\frac{1}{288}\, T_4$, 
the equation
$\overset{\green{\bf 6001}}{\,=\,}$ above
would be a nontrivial linear dependence relation
between the transitivity parameters $T_1, T_2, T_3, T_4$,
which is forbidden.
\endproof

From $\eqFG$:
\[
0
\overset{\green{\bf 6001}}{\,\,=\,\,}
-\,\tfrac{1}{720}\,F_{6,0,0,1}\,a_{1,1}^2
+
\tfrac{1}{720}\,G_{6,0,0,1}\,a_{1,1}^4,
\]
we see that we can normalize:
\[
G_{6,0,0,1}
\,:=\,
1
\ \ \ \ \ \ \ \ \ \ \ \ \ \ \ \ \ \ \ \
\text{or}
\ \ \ \ \ \ \ \ \ \ \ \ \ \ \ \ \ \ \ \
G_{6,0,0,1}
\,:=\,
-\,1,
\]
and the same about $F_{6,0,0,1}$.
Stabilization of this last normalization
requires $a_{1, 1} := 1$, and
if there exists any homogeneous model,
it can only be simply transitive.

As a first case, 
put $F_{6,0,0,1} := 1$ everywhere, solve from
$\overset{\green{\bf 6001}}{\,=\,}$:
\[
A_{1,1}
\,:=\,
\Big(
-\tfrac{3}{2}\,F_{5,1,0,0}
-
\tfrac{1}{2}\,F_{7,0,0,1}\,T_1
\Big)\,
T_1
-
T_2
+
\tfrac{5}{4}\,T_4,
\]
and replace this value of $A_{1,1}$ everywhere.
Then $\overset{\green{\bf 5100}}{\,=\,}$ becomes:
\[
0
\overset{\green{\bf 5100}}{\,\,=\,\,}
\Big(
\tfrac{1}{120}\,F_{6,1,0,0}
-
\tfrac{1}{80}\,F_{5,1,0,0}\,F_{7,0,0,1}
-
\tfrac{3}{80}\,F_{5,1,0,0}^2
\Big)\,T_1
+
\tfrac{1}{120}\,F_{5,1,0,0}\,T_2
+
\tfrac{1}{32}\,F_{5,1,0,0}\,T_4,
\]
whence necessarily:
\[
F_{5,1,0,0}
\,=\,
0
\ \ \ \ \ \ \ \ \ \ \ \ \ \ \ \ \ \ \ \
\text{and then:}
\ \ \ \ \ \ \ \ \ \ \ \ \ \ \ \ \ \ \ \
F_{6,1,0,0}
\,=\,0.
\]
This necessary vanishing $F_{5,1,0,0} = 0$ {\em a posteriori}
explains why we did not open a branch {\em supra}.

Therefore, put $F_{5,1,0,0} := 0$ and $F_{6,1,0,0} := 0$ everywhere.
Then $\overset{\green{\bf 7000}}{\,=\,}$ becomes:
\[
\aligned
0
&
\overset{\green{\bf 7000}}{\,\,=\,\,}
\Big(
\tfrac{1}{5040}\,F_{8,0,0,0}
-
\tfrac{1}{3600}\,F_{7,0,1,0}
-
\tfrac{1}{2016}\,F_{7,0,0,0}\,F_{7,0,0,1}
-
\tfrac{1}{1900}
\Big)\,T_1
\\
&
\ \ \ \ \
+
\Big(
\tfrac{1}{5040}\,F_{7,1,0,0}
-
\tfrac{1}{840}\,F_{7,0,0,0}
\Big)\,T_2
+
\Big(
\tfrac{1}{5040}\,F_{7,0,1,0}
+
\tfrac{1}{5400}
\Big)\,T_3
+
\Big(
\tfrac{1}{5040}\,F_{7,0,0,1}
+
\tfrac{5}{4032}\,F_{7,0,0,0}
\Big)\,T_4.
\endaligned
\]
It follows:
\[
F_{7,0,1,0}
\,=\,
-\,\tfrac{5040}{5400}
\,=\,
-\,\tfrac{14}{15},
\]
whence:
\[
0
\overset{\green{\bf 7001}}{\,\,=\,\,}
\Big(
\tfrac{1}{5040}\,F_{8,0,0,1}
-
\tfrac{1}{3360}\,F_{7,0,0,1}^2
\Big)\,T_1
+
\tfrac{1}{1344}\,F_{7,0,0,1}\,T_4,
\]
so that:
\[
F_{7,0,0,1}
\,=\,
0,
\ \ \ \ \ \ \ \ \ \ \ \ \ \ \ \ \ \ \ \
F_{8,0,0,1}
\,=\,
0.
\]

Then:
\[
\aligned
0
\overset{\green{\bf 7000}}{\,\,=\,\,}
\Big(
\tfrac{1}{5040}\,F_{8,0,0,0}
+
\tfrac{1}{6750}
\Big)\,T_1
+
\Big(
\tfrac{1}{5040}\,F_{7,1,0,0}
-
\tfrac{1}{840}\,F_{7,0,0,0}
\Big)\,T_2
+
\tfrac{5}{4032}\,F_{7,0,0,0}
\,T_4.
\endaligned
\]
Thus:
\[
\aligned
F_{7,0,0,0}
&
\,=\,
0,
&
\ \ \ \ \ \ \ \ \ 
F_{7,1,0,0}
&
\,=\,
0,
&
\ \ \ \ \ \ \ \ \ 
F_{8,0,0,0}
&
\,=\,
-\,
\tfrac{56}{75},
\\
F_{7,1,0,0}
&
\,=\,
0,
&
\ \ \ \ \ \ \ \ \ 
F_{8,0,1,0}
&
\,=\,
0,
&
\ \ \ \ \ \ \ \ \ 
F_{8,1,0,0}
&
\,=\,
-\,\tfrac{392}{75},
\ \ \ \ \ \ \ \ \ 
F_{9,0,0,0}
\,=\,
0.
\endaligned
\]

The second case $F_{6,0,0,1} = - 1$ is treated similarly.

\begin{Theorem}
Among constant Hessian rank 1 hypersurfaces
$H^4 \subset \R^5$, there are only two 
affinely homogeneous models, 
of equations depending on some sign choices
$\pm$ or $\mp$:
\[
\aligned
u
&
\,=\,
\ \
\tfrac{x^2}{2}
\\
&
\ \ \ \ \
+
\tfrac{x^2y}{2}
\\
&
\ \ \ \ \
+
\tfrac{x^3z}{6}
+
\tfrac{x^2y^2}{2}
\\
&
\ \ \ \ \
+
\tfrac{x^4w}{24}
+
\tfrac{x^3yz}{2}
+
\tfrac{x^2y^3}{2}
\\
&
\ \ \ \ \
+
\tfrac{x^4yw}{6}
+
\tfrac{x^4z^2}{8}
+
x^3y^2z
+
\tfrac{x^2y^4}{2}
\\
&
\ \ \ \ \
\pm\,\tfrac{x^6w}{720}
+
\tfrac{1}{12}\,x^5zw
+
\tfrac{5}{12}\,x^4y^2w
+
\tfrac{5}{8}\,x^4yz^2
+
\tfrac{5}{3}\,x^3y^3z
+
\tfrac{1}{2}\,x^2y^5
\\
&
\ \ \ \ \
\mp\,\tfrac{x^8}{54000}
-
\tfrac{x^7z}{5400}
\pm\,\tfrac{x^6yw}{120}
+
\tfrac{x^6w^2}{72}
+
\tfrac{x^5yzw}{2}
+
\tfrac{x^5z^3}{8}
+
\tfrac{5}{6}\,x^4y^3w
+
\tfrac{15}{8}\,x^4y^2z^2
+
\tfrac{5}{2}\,x^3y^4z
+
\tfrac{x^2y^6}{2}
\\
&
\ \ \ \ \
+
{\rm O}_{x,y,z,w}(9),
\endaligned
\]
with 4-dimensional affine symmetry algebra generated by:
\[
\aligned
e_1
&
\,:=\,
\big(
1-y\pm\tfrac{1}{5}u
\big)\,\partial_x
+
\big(
\mp\tfrac{1}{5}x-z
\big)\,\partial_y
+
\big(
-w-\tfrac{4}{75}u
\big)\,\partial_z
+
\big(
\tfrac{8}{75}x\pm\tfrac{2}{5}z
\big)\,\partial_w
+
x\,\partial_u,
\\
e_2
&
\,:=\,
-\,x\partial_w
+
(1-y)\,\partial_y
-
z\,\partial_z
-
w\,\partial_w
-
u\,\partial_u,
\\
e_3
&
\,:=\,
\tfrac{2}{3}\,u\,\partial_x
-
x\,\partial_y
+
\big(
1-y\mp\tfrac{1}{15}u
\big)\,\partial_z
+
\big(
\pm\tfrac{2}{15}x-\tfrac{2}{3}z
\big)\,\partial_w
\\
e_4
&
\,:=\,
\pm\,\tfrac{5}{4}\,x\,\partial_x
+
\tfrac{1}{2}\,u\,\partial_y
+
\big(
-x+\tfrac{5}{4}\,z
\big)\,\partial_z
+
\big(
1-y\mp\tfrac{5}{2}w\mp\tfrac{1}{15}u
\big)\,\partial_w
\pm
\tfrac{5}{2}\,u\,\partial_u,
\endaligned
\]
sharing the Lie brackets:
\reqnomode\usetagform{EngelLie}
\begin{align}
{}
[e_1,e_2]
\,=\,
0,
\ \ \ \ \ \ \ \ \ \ \ 
[e_1,e_3]
\,=\,
\mp\,\tfrac{4}{15}\,e_4,
\ \ \ \ \ \ \ \ \ \ \ 
[e_1,e_4]
&
\,=\,
\pm\,\tfrac{5}{4}\,e_1,
\notag
\\
{}
[e_2,e_3]
\,=\,
0,
\ \ \ \ \ \ \ \ \ \ \ 
[e_2,e_4]
&
\,=\,
0,
\notag
\\
{}
[e_3,e_4]
&
\,=\,
\mp\,\tfrac{5}{4}\,e_3.
\tag{\qed}
\end{align}
\end{Theorem}



\vfill\end{document}